\documentclass[journal]{IEEEtran}

\usepackage{color}

\ifCLASSINFOpdf
   \usepackage[pdftex]{graphicx}
\else
\fi
\usepackage{amsmath}
\usepackage{amsfonts,amssymb,amsthm}
\usepackage{algorithmic}
\usepackage{hyperref}
\usepackage{amssymb}
\usepackage{tabularx}
\usepackage{multirow}
\usepackage{booktabs}
\usepackage{pifont}

\usepackage{arydshln}
\usepackage{color}
\usepackage{xcolor}

\newsavebox{\ieeealgbox}
\newenvironment{boxedalgorithmic}
  {\begin{lrbox}{\ieeealgbox}
   \begin{minipage}{\dimexpr\columnwidth-2\fboxsep-2\fboxrule}
   \begin{algorithmic}}
  {\end{algorithmic}
   \end{minipage}
   \end{lrbox}\noindent\fbox{\usebox{\ieeealgbox}}}

\newtheorem{theorem}{Theorem}[section]

\newtheorem{definition}[theorem]{Definition}

\newcommand{\modNT}[1]{\textcolor{black}{#1}}
\newcommand{\modRG}[1]{\textcolor{black}{#1}}
\newcommand{\modLLM}[1]{\textcolor{black}{#1}}
\newcommand{\modfinal}[1]{\textcolor{black}{#1}}

\begin{document}
%
\title{\modLLM{Approximate fast graph Fourier transforms\\ via multi-layer sparse approximations}}
%
%
%

\author{Luc~Le Magoarou,
        R\'emi~Gribonval,~\IEEEmembership{Fellow,~IEEE},
        Nicolas Tremblay
\thanks{Luc~Le Magoarou (\href{mailto:luc.le-magoarou@inria.fr}{\texttt{luc.le-magoarou@inria.fr}}) and R\'emi Gribonval (\href{mailto:remi.gribonval@inria.fr}{\texttt{remi.gribonval@inria.fr}}) are both with Inria, Rennes, France, PANAMA team. Nicolas Tremblay is with CNRS, GIPSA-lab, Grenoble, France. This work was supported in part by the European Research Council, PLEASE project (ERC-StG- 2011-277906). 
}
}

%
%

\markboth{}%
{Shell \MakeLowercase{\textit{et al.}}: Bare Demo of IEEEtran.cls for Journals}
%



\maketitle

\begin{abstract}
The Fast Fourier Transform (FFT) is an algorithm of paramount importance in signal processing as it allows to apply the Fourier transform in $\mathcal{O}(n \log n )$ instead of $\mathcal{O}(n^2)$ arithmetic operations. Graph Signal Processing (GSP) is a recent research domain that generalizes classical signal processing tools, such as the Fourier transform, to situations where the signal domain is given by any arbitrary graph instead of a regular grid. Today, there is no method to rapidly apply graph Fourier transforms. We propose in this paper a method to obtain approximate graph Fourier transforms that can be applied rapidly  \modRG{and stored efficiently}. It is based on a greedy approximate diagonalization of the graph Laplacian matrix, carried out using a modified version of the famous Jacobi eigenvalues algorithm. The method is described and analyzed \modLLM{in detail}, and then applied to both synthetic and real graphs, showing its potential.
\end{abstract}


%
\IEEEpeerreviewmaketitle

\section{Introduction}
%
%
%
%
\IEEEPARstart{G}{raphs} are powerful mathematical objects used to model pairwise relationships between elements of a set. Graph theory has been extensively developed since the eighteenth century, and has found a variety of applications, ranging from biology to computer science or linguistics \cite{Bondy1976}.

Recently, methods have been developed to analyze and process signals defined over the vertices of a graph \cite{Shuman2013,Sandryhaila2013}, instead of over a regular grid, as is classically assumed in discrete signal processing. The starting point of graph signal processing is to define a Fourier transform, via an analogy with classical signal processing. Depending on the preferred analogy, there exists several definitions of graph Fourier transforms \cite{Shuman2013,Sandryhaila2013}. Following \cite{Shuman2013}, we choose in this paper to define the graph Fourier basis as the eigenvector  matrix of the graph Laplacian operator $\mathbf{L}$ (whose precise definition is given in \modLLM{Section}~\ref{sec:notations}). 
In a graph with $n$ vertices, $\mathbf{L} \in \mathbb{R}^{n \times n}$ and 
\begin{equation}
\mathbf{L} = \mathbf{U}\boldsymbol{\Lambda}\mathbf{U}^T,
\label{eq:defGFT}
\end{equation} 
where $\mathbf{U} \in \mathbb{R}^{n \times n}$ is an orthogonal matrix whose columns are the graph Fourier modes and $\boldsymbol{\Lambda} \in \mathbb{R}^{n \times n}$ is a \modLLM{non-negative} diagonal matrix whose diagonal entries correspond to the 
squared graph frequencies. 
 
 The graph Fourier matrix $\mathbf{U}$ being non-sparse in general, applying it costs $\mathcal{O}(n^2)$ arithmetic operations. In the classical signal processing case, the well-known Fast Fourier Transform (FFT) \cite{CooleyTukey1965} allows to apply the Fourier transform in only $\mathcal{O}(n \log n)$ arithmetic operations. The FFT is a fast linear algorithm \cite{Morgenstern1975}, which implies that the classical Fourier matrix can be factorized into sparse factors, as discussed in \cite{Lemagoarou2016a}. Given that classical signal processing is equivalent to graph signal processing on the  ring graph, it is natural to wonder if this kind of factorization can be generalized to other graphs.

We proposed in a previous work a method to obtain computationally efficient\footnote{\modLLM{An $m \times n$ matrix is ``efficient'' if it is associated with a linear algorithm involving strictly less than  $mn$ scalar multiplications.}} approximations of matrices, based on multi-layer sparse factorizations \cite{Lemagoarou2016a}.  
 The method amounts to approximate a matrix of interest $\mathbf{A} \in \mathbb{R}^{p \times q}$ with a product of sparse matrices, as 
\begin{equation}
\mathbf{A} \approx \mathbf{S}_J\dots\mathbf{S}_1,
\label{eq:spop}
\end{equation} 
where the matrices $\mathbf{S}_1,\dots,\mathbf{S}_J$ are sparse, allowing for cheap storage and multiplication. 
We applied this method to the graph Fourier matrix $\mathbf{U}$ of various graphs in \cite{Lemagoarou2016}, in order to get \modLLM{approximate Fast Graph Fourier Transforms
(FGFT)}
The method showed good potential. However in the context of graph signal processing, this approach suffers from at least two limitations:
 \begin{enumerate}
 \item[\bf(L1)] It requires a full diagonalization of the graph Laplacian $\mathbf{L}$ before it can be applied. Indeed, the method takes the graph Fourier matrix $\mathbf{U}$ as input. Performing this diagonalization costs $\mathcal{O}(n^3)$ arithmetic operations, which is prohibitive if $n$ is \modNT{large}. 
 \item[\bf(L2)] It provides non-orthogonal approximate \modNT{FGFTs}. 
 Indeed, the details of the method make it difficult to get sparse \emph{and} orthogonal factors $\mathbf{S}_1,\dots,\mathbf{S}_J$. This leads to approximate graph Fourier transforms that are not easily invertible, which can be a problem for applications where signal reconstruction is needed.
 \end{enumerate}

We propose in this paper another method, that does not suffer from these limitations, to obtain approximate \modLLM{FGFT}s. In order to overcome {\bf(L1)}, we consider directly the Laplacian matrix $\mathbf{L}$ as input. To overcome {\bf(L2)}, we look for factors $\mathbf{S}_1,\dots,\mathbf{S}_J$ constrained to be in a set of sparse and orthogonal matrices built with Givens rotations \cite{Givens1958} (as explained \modLLM{in detail} in section~\ref{sec:opt_fram}). The proposed method amounts to an approximate diagonalization of the Laplacian matrix $\mathbf{L}$, as
\begin{equation}
\mathbf{L} \approx \mathbf{S}_1\dots\mathbf{S}_J\hat{\boldsymbol{\Lambda}}\mathbf{S}_J^T\dots\mathbf{S}_1^T,
\label{eq:defGFFT}
\end{equation} 
where the matrices $\mathbf{S}_1,\dots,\mathbf{S}_J$ are both sparse \emph{and} orthogonal. The product $\hat{\mathbf{U}} = \mathbf{S}_1\dots\mathbf{S}_J$ constitutes an efficient approximate graph Fourier matrix and $\hat{\boldsymbol{\Lambda}}$ is a diagonal matrix whose diagonal entries are approximations of the squared graph frequencies.\\

\begin{figure}[tb]
\centering
\includegraphics[width=\columnwidth]{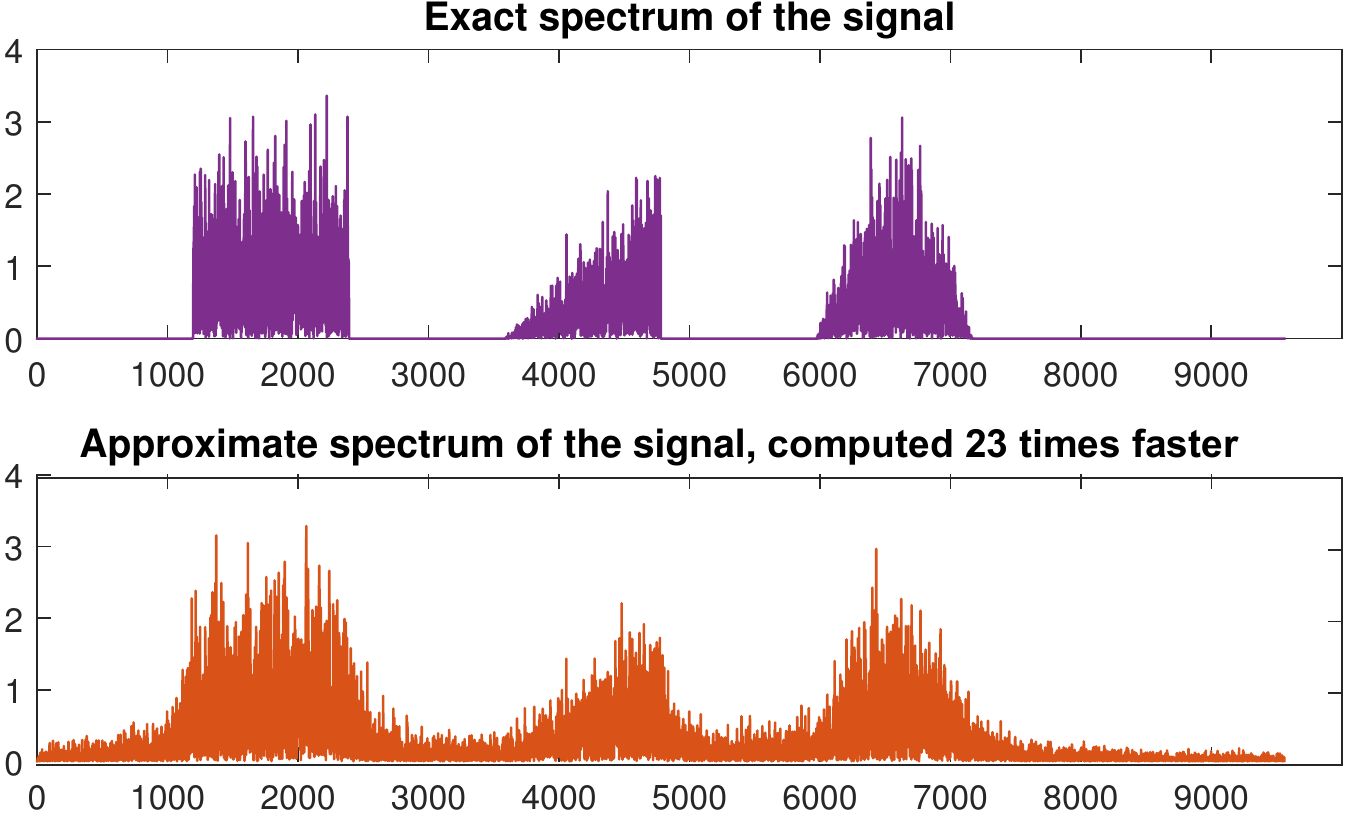}
\caption{Example of an approximate spectrum computed with our method on a real-world sensor graph. Details in Section~\ref{sec:sensor}.}
\label{fig:spectrumexample}
\end{figure}

\noindent {\bf Contributions.}
Given a graph Laplacian matrix, the main objective of this paper is to find an approximate graph Fourier matrix $\hat{\mathbf{U}}$ that both i) approximately diagonalizes the Laplacian and ii) is computationally efficient. 
The proposed method is a greedy strategy that amounts to \modLLM{truncating and slightly modifying} the famous Jacobi eigenvalue algorithm \cite{Jacobi1846}. \modLLM{Note that it could in principle be applied to any symmetric matrix (covariance matrix, normalized Laplacian, etc.), but the focus of the present paper is on the graph combinatorial Laplacian. 
Indeed, while a general symmetric matrix has no reason to be well approximated with a limited number of Givens factors, graph Laplacians are shown empirically to have this interesting property \modfinal{for several popular graph families}.} The proposed method is compared experimentally with the direct approximation method of \cite{Lemagoarou2016} on various graphs, showing that it is much faster while obtaining good results. In fact, we obtain  approximate \modLLM{FGFTs} of complexity $\mathcal{O}(n\log n)$ exhibiting constant error for growing $n$. \modLLM{Moreover, we briefly investigate the links between the graph structure and the quality of diagonalization.} Also, a discussion on possible use cases of the method is undertaken and further experimental validations on both synthetic and real-world sensor networks are conducted, showing an interesting compromise between computational complexity and accuracy. \modLLM{Figure}~\ref{fig:spectrumexample} illustrates one such experiment further developed in Section~\ref{sec:sensor}.  
Finally, the method is applied to graph signal filtering, and compared to the usual method based on polynomial approximations. We show that when the filter function is not well approximated by a polynomial (such as low-pass filters with a steep cut-off), our method performs at least as \modLLM{well}, and sometimes outperforms the polynomial approximation method.\\ 

The remainder of the paper is organized as follows. The problem is formulated in \modLLM{Section~}\ref{sec:prob_form}, and an approach to solve it is given in \modLLM{Section~}\ref{sec:opt_fram}. \modLLM{The proposed method is then evaluated, compared to the direct approximation method and its links with the graph structure are briefly investigated in \modLLM{Section~}\ref{sec:expval}}. Then, an application to sensor networks is proposed in \modLLM{Section~}\ref{sec:sensor}. Finally, an application to graph signal filtering is presented in \modLLM{Section~}\ref{sec:filtering}.

\section{Problem formulation}
\label{sec:prob_form}
In this section, we give a concrete formulation of the main problem \modNT{addressed in this} paper. We first set up the notations and conventions used, before presenting the objective in detail. We then discuss the advantages in terms of computational complexity expected from the method and end with a presentation of the related work.

\subsection{Notations and conventions}
\label{sec:notations}

{\noindent \bf General notations.} Matrices are denoted by bold upper-case letters: $\mathbf{A}$; vectors  by bold lower-case letters: $\mathbf{a}$; the $i$th column of a matrix $\mathbf{A}$ by: $\mathbf{a}_i$; its entry located at the $i$th row and $j$th column by: $a_{ij}$. 
Sets are denoted by calligraphic symbols: $\mathcal{A}$, and we denote by $\delta_{\mathcal{A}}$ the \modfinal{characteristic} function of the set $\mathcal{A}$ in the optimization sense ($\delta_{\mathcal{A}}(x) = 0$ if $x \in \mathcal{A}$, $\delta_{\mathcal{A}}(x) = +\infty$ otherwise). The standard vectorization operator is denoted $\text{vec}(\cdot)$. The $\ell_0$-pseudonorm is denoted $\left\Vert\cdot\right\Vert_0$ (it counts the number of non-zero entries), $\left\Vert\cdot\right\Vert_F$ denotes the Frobenius norm, and $\left\Vert\cdot\right\Vert_{2}$ the spectral norm. By abuse of notations, $\|\mathbf{A}\|_{0} = \|\text{vec}(\mathbf{A})\|_{0}$. $\mathbf{Id}$ denotes the identity matrix. \\

{\noindent \bf Graph Laplacian.} We consider in this paper undirected weighted graphs, denoted $\mathcal{G} \triangleq \{\mathcal{V},\mathcal{E},\mathbf{W}\}$, where $\mathcal{V}$ represents the set of vertices (otherwise called nodes), $\mathcal{E} \subset \mathcal{V} \times \mathcal{V}$ is the set of edges, and $\mathbf{W}$ is the weighted adjacency matrix of the graph. We denote $n\triangleq |\mathcal{V}|$ the total number of vertices and the adjacency matrix $\mathbf{W} \in \mathbb{R}^{n \times n}$ is symmetric and such that $w_{ij} = w_{ji}$ is non-zero only if $(i,j)\in \mathcal{E}$ and represents the strength of the connection between nodes $i$ and $j$. We define the degree matrix $\mathbf{D} \in \mathbb{R}^{n \times n}$ as a diagonal matrix with $\forall i$, $d_{ii} \triangleq \sum_{j=1}^n w_{ij}$, and the combinatorial Laplacian matrix $\mathbf{L}\triangleq \mathbf{D} - \mathbf{W}$ (we only consider this type of Laplacian matrix in this paper, and hereafter simply call it the Laplacian).

\subsection{Objective}
Our goal is to approximately diagonalize the Laplacian $\mathbf{L}$ with an efficient approximate eigenvector matrix $\hat{\mathbf{U}} = \mathbf{S}_1\dots \mathbf{S}_J$, where the factors $\mathbf{S}_1,\dots,\mathbf{S}_J \in \mathbb{R}^{n \times n}$ are sparse and orthogonal. Using the Frobenius norm to measure the quality of approximation, this objective can be stated as the following optimization problem:
\begin{equation}\tag{DP}
\begin{array}{cl}
\underset{\hat{\boldsymbol{\Lambda}}, \mathbf{S}_1,\dots,\mathbf{S}_J}{\text{minimize}}& \quad \left\Vert \mathbf{L} - \mathbf{S}_1\dots \mathbf{S}_J \hat{\boldsymbol{\Lambda}} \mathbf{S}_J^T\dots \mathbf{S}_1^T \right\Vert_F^2  
\\
& \quad + \sum_{j=1}^J \delta_{\mathcal{S}}(\mathbf{S}_j) + \delta_{\mathcal{D}}(\hat{\boldsymbol{\Lambda}}),
\end{array}
\label{eq:diago}
\end{equation}
where $\mathcal{D}$ is a set of diagonal matrices and $\mathcal{S}$ is a set of sparse and orthogonal matrices. The set of sparse matrices $\{\mathbf{A} \in \mathbb{R}^{n \times n}, \left\Vert \mathbf{A} \right\Vert_0 \leq s\}$ and the set of orthogonal matrices $\{\mathbf{A} \in \mathbb{R}^{n \times n}, \mathbf{A}^T\mathbf{A} = \mathbf{Id}\}$ are both easy to project onto, but their intersection is not. This is the reason why approaches using projected gradient descent similar to what is proposed in \cite{Lemagoarou2016a} cannot be used to solve this problem. Instead, we propose a greedy strategy using Givens rotations that is explained \modLLM{in detail} in \modLLM{Section~}\ref{sec:opt_fram}.

\subsection{Relative Complexity Gain}

An efficient approximate Fourier matrix $\hat{\mathbf{U}} = \mathbf{S}_1\dots \mathbf{S}_J$ is beneficial in terms of complexity provided its \emph{relative complexity} with respect to the exact Fourier matrix $\mathbf{U}$ is small. In order to quantify this statement, we introduce a measure of the relative complexity in the following definition.

\begin{definition}
The Relative Complexity Gain (abbreviated RCG) is the ratio between the number of non-zero entries in the Fourier matrix $\mathbf{U}$ and the total number of non-zero entries in the multi-layer sparse factorization of the efficient approximate Fourier matrix $\hat{\mathbf{U}} = \mathbf{S}_1\dots \mathbf{S}_J$ : 
\begin{equation}
\text{RCG} \triangleq \frac{\left\Vert \mathbf{U} \right\Vert_0}{\sum_{j=1}^J\left\Vert \mathbf{S}_j \right\Vert_0} .
\label{RC} 
\end{equation} 
\end{definition}

The relative complexity gain is a good measure of the \emph{theoretical} benefit of using the efficient approximation $\hat{\mathbf{U}} = \mathbf{S}_1\dots \mathbf{S}_J$ in place of the exact Fourier matrix $\mathbf{U}$, in several \modLLM{aspects} of its manipulation, such as storage and multiplication (see \cite{Lemagoarou2016a} for a detailed explanation on this). For example, having  $\text{RCG} = 10$ means that the multiplication of any vector $\mathbf{x} \in \mathbb{R}^n$ by the efficient approximation or its transpose will require ten times less scalar multiplications than a multiplication by $\mathbf{U}$. \modRG{It also means that storing the approximate FGFT will require roughly ten times less memory than the original Fourier transform, as will be demonstrated on the DC road graph which Fourier matrix requires 700 MB while its  approximate FGFT only requires 60MB.} However, this does not necessarily mean that the multiplication will be ten times faster \emph{in practice}, as we will see in the experimental part of this paper (Sections~\ref{sec:expval} to~\ref{sec:filtering}.)

\subsection{Related work}

Computationally efficient approximations of operators have been studied quite extensively in various domains over the last years. Applications of such methods range from dictionary learning \cite{Rubinstein2010a, Chabiron2014, Sulam2016, Chabiron2016}, to Non-negative Matrix Factorization (NMF) \cite{Lyu2013} or compression of integral operators arising from discretized continuous kernels \cite{Beylkin1991,Rokhlin1985,Hackbusch1999,Candes2007}.
For a brief description of these various methods in a unified way, see \cite{Lemagoarou2016a}.

We focus here on efficient \emph{orthogonal} approximation of operators. In this context, the idea to use a multi-layer product of sparse and orthogonal matrices (otherwise called elementary rotations) in order to make complexity savings 
has been investigated in several approaches over the last few years.  Very recently, it has been proposed \cite{Rusu2016} to learn an efficient orthogonal dictionary being the product of a low number of Householder reflectors \cite{Householder1958} (another type of sparse and orthogonal matrices, different from Givens rotations).

Several approaches have also considered diagonalizing symmetric matrices with a product of sparse and orthogonal matrices.
The first one comes from statistics \cite{Lee2008}, and introduces a method to design an efficient orthogonal transform to decorrelate some data, based on an empirical similarity matrix that is to be approximately diagonalized with a low number of Givens rotations \cite{Givens1958} whose supports are organized as a tree. Another approach \cite{Cao2011} aims at estimating the covariance matrix of some data, under the constraint that its eigenvector matrix is the product of a low number of Givens rotations, with a maximum likelihood formulation \modLLM{involving determinants}. \modLLM{In the context of video compression, it has also been proposed \cite{Chen2011,Chen2012} to optimize the coding gain by an approximate diagonalization of the covariance matrix with Givens rotations, with a criterion involving logarithms}.

 More related to graphs and to the method we propose in this paper, two recent approaches called Multiresolution Matrix Factorization (MMF) \cite{Kondor2014} and parallel Multiresolution Matrix Factorization (pMMF) \cite{Teneva2016} propose to design a wavelet transform on graph via a greedy factorization of the graph Laplacian with the help of general elementary rotations being constrained to comply with certain multiresolution axioms. For the pMMF approach, the factorization algorithm is parallelized using a clustering algorithm that cuts down the matrix to factorize into several pieces, and is applied several times along the factorization algorithm. MMF and pMMF have been used for various tasks such as matrix compression and preconditioning.

\modLLM{A complementary approach to the design of fast graph Fourier transforms is presented in \cite{Sandryhaila2014}, where certain structured graphs are viewed as the product of two (or more) subgraphs (the whole graph being called a \emph{product graph}). The product structure can be used to reduce the Fourier transform complexity. The two approaches are perfectly compatible, since the Fourier transform on the product graph can be further accelerated by computing approximate FGFTs for the Laplacian matrix of one (or more) of the subgraphs. For example if a graph represents the variations along time of a signal defined on a sensor network, it can be seen as the product of a sensor graph and a time series. The Laplacian matrix of the sensor graph can be approximately diagonalized, while the product structure can still be leveraged to decouple it from the time variations and accelerate the overall Fourier transform.}

The approach we propose in this paper amounts to an approximate diagonalization with a product of sparse and orthogonal matrices. However it differs from \cite{Lee2008,Cao2011} since it departs from the probabilistic formulation that induces different update rules, and is applied to graph Laplacians rather than covariances or similarity matrices. On the other hand, it is closer to the MMF approaches \cite{Kondor2014,Teneva2016}, but the main objectives and targeted applications differ. Indeed, we are only interested here in having an efficient linear transformation that approximately diagonalizes the graph Laplacian, so that we do not impose constraints during the approximate diagonalization to achieve some form of multiresolution. In that sense, it can be seen as simpler.

In summary, and as will be made clearer in the next section, the method we propose here can be seen as an approximate diagonalization using \emph{a modified (in particular, truncated) 
version of the Jacobi eigenvalues algorithm} \cite{Jacobi1846,Golub2000}, where the truncation allows for both a faster computation of the approximate Fourier transform and a faster multiplication by this approximate transform, as demonstrated numerically in Section~\ref{sec:expval}. \modLLM{Actually, the novelty and genuine contributions of this paper primarily consist in demonstrating numerically that the described algorithms can be used to get approximate \modLLM{FGFTs}. Indeed, the problem of \modLLM{computing rapidly the graph Fourier transform} has not yet been studied extensively, despite its great importance.}


\section{Optimization framework}
\label{sec:opt_fram}

In this section, we describe the strategy we propose to solve the optimization problem~\eqref{eq:diago} introduced in \modLLM{Section~}\ref{sec:prob_form}.

\subsection{Constraint sets}
We impose that each sparse factor $\mathbf{S}_j$ belong to the set $\mathcal{R}_G$ of Givens rotations \cite{Givens1958} ($\mathcal{S} = \mathcal{R}_G$). An $n$-dimensional Givens rotation is a linear transformation that does not act on $n-2$ coordinates and rotates the two remaining by an angle $\theta\in [0;2\pi[$. Noting $p$ and $q$ the two rotated coordinates, Givens rotations thus correspond to matrices of the following form,

%

\begin{figure}[!h]
\centering
\includegraphics[width=0.9\columnwidth]{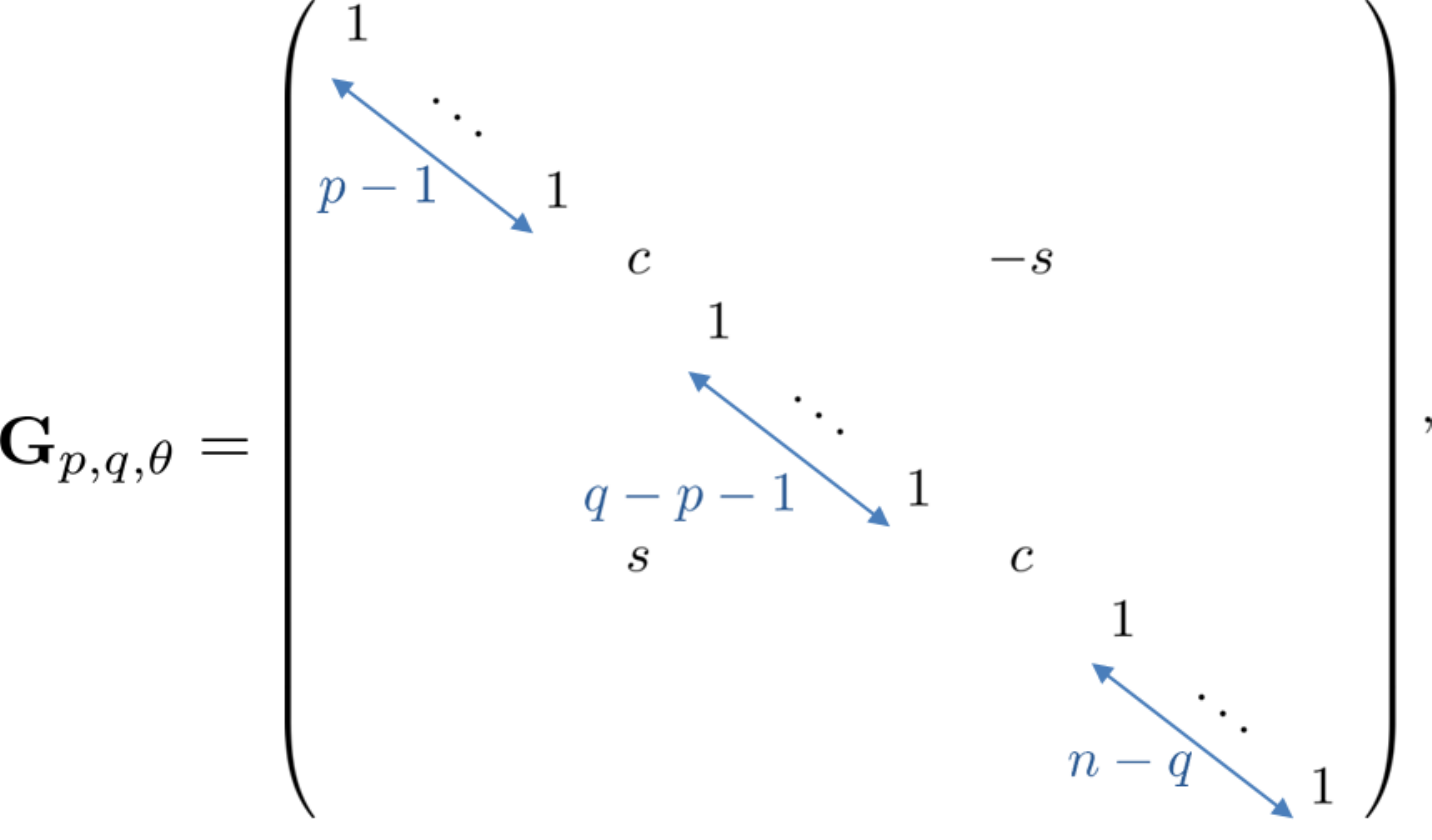}
\end{figure}
\noindent where $c\triangleq\cos(\theta)$ and $s \triangleq \sin(\theta)$. A Givens rotation only depends on three parameters: the two coordinates $p$ and $q$ and the rotation angle $\theta$, hence the notation $\mathbf{G}_{p,q,\theta}$. A Givens rotation can be considered as a two-dimensional linear transformation, so that in terms of computational complexity, and for the Relative Complexity Gain (RCG) calculation, we consider that $\left\Vert \mathbf{G}_{p,q,\theta} \right\Vert_0 =4$. Regarding the matrix $\hat{\boldsymbol{\Lambda}}$ containing the estimated squared graph frequencies, we do not impose any constraint on it except that it should be diagonal. This corresponds to taking $\mathcal{D}$ as the set of all diagonal matrices.
 
\subsection{Truncated Jacobi approach} 
To approximately solve problem~\eqref{eq:diago}, one can rely on a truncated Jacobi eigenvalues algorithm \cite{Jacobi1846,Golub2000}. The Jacobi algorithm is an iterative procedure where at each step one seeks the Givens rotation that reduces the most the cost function. At the first step, this means setting $\mathbf{S}_1$ and $\hat{\boldsymbol{\Lambda}}$ as follows:
\begin{equation*}
(\mathbf{S}_1,\hat{\boldsymbol{\Lambda}}) \leftarrow \underset{\mathbf{D} \in \mathcal{D}, \mathbf{S} \in \mathcal{R}_G}{\text{argmin}} \quad \left\Vert \mathbf{L} - \mathbf{S} \mathbf{D}  \mathbf{S}^T \right\Vert_F^2,
\end{equation*}
which can be reformulated using the fact that the Frobenius norm is invariant under orthogonal transformations as
\begin{equation*}
(\mathbf{S}_1,\hat{\boldsymbol{\Lambda}}) \leftarrow \underset{\mathbf{D} \in \mathcal{D}, \mathbf{S} \in \mathcal{R}_G}{\text{argmin}} \quad \left\Vert \mathbf{S}^T\mathbf{L} \mathbf{S}- \mathbf{D}   \right\Vert_F^2.
\end{equation*}
The set $\mathcal{D}$ being the set of all diagonal matrices, the optimal estimated squared graph frequencies factor is simply $\hat{\boldsymbol{\Lambda}} = \text{diag}(\mathbf{S}^T\mathbf{L} \mathbf{S})$. This allows to rule out this factor of the problem and to reformulate it as follows:
\begin{equation*}
\mathbf{S}_1 \leftarrow\underset{ \mathbf{S} \in \mathcal{R}_G}{\text{argmin}} \quad \left\Vert \mathbf{S}^T\mathbf{L} \mathbf{S}   \right\Vert_{\text{offdiag}}^2,
\end{equation*} 
where $\left\Vert \mathbf{A}\right\Vert_{\text{offdiag}}^2$  is simply the sum of the squared off-diagonal entries of $\mathbf{A}$. Once the factor $\mathbf{S}_1$ is set this way, and introducing the notation $\mathbf{L}_2 \triangleq \mathbf{S}_1^T\mathbf{L} \mathbf{S}_1$, the next step of the strategy is to choose  $\mathbf{S}_2 \leftarrow\underset{ \mathbf{S} \in \mathcal{R}_G}{\text{argmin}} \quad \left\Vert \mathbf{S}^T\mathbf{L}_2 \mathbf{S}   \right\Vert_{\text{offdiag}}^2$, and so on until the $J$th and last step. The algorithm thus amounts to solve a sequence of $J$ very similar subproblems of the form 
\begin{equation}\tag{SP}
\underset{ \mathbf{S} \in \mathcal{R}_G}{\text{minimize}} \quad \left\Vert \mathbf{S}^T\mathbf{L}_j \mathbf{S}   \right\Vert_{\text{offdiag}}^2,
\label{eq:subprob}
\end{equation}
with $\mathbf{L}_j \triangleq \mathbf{S}_{j-1}^T \mathbf{L}_{j-1}\mathbf{S}_{j-1}$. This is summarized by the algorithm of \modLLM{Figure}~\ref{algo:givens_sum}. Compared to the traditional Jacobi eigenvalues algorithm \cite{Jacobi1846,Golub2000}, where new Givens rotations are chosen until a certain accuracy is attained, the main difference is that by prescribing the number $J$ of Givens rotations we can adjust the trade-off between accuracy and computational efficiency of the product $\mathbf{S}_1\dots\mathbf{S}_J$.

\modfinal{
Note that in order for the approximate FGFT $\hat{\mathbf{U}} = \mathbf{S}_1\dots\mathbf{S}_J$ to make sense in a graph signal processing context, we reorder its columns according to the estimated eigenvalues (line $7$ of the algorithm). This way, the first columns of $\hat{\mathbf{U}} = \mathbf{S}_1\dots\mathbf{S}_J$ ``physically'' correspond to low frequencies, and its last columns to high frequencies.
}

\begin{figure}[tbp]
\centering
\begin{boxedalgorithmic}[1] 
\REQUIRE{matrix $\mathbf{L}$, number $J$ of Givens rotations.}\STATE $\mathbf{L}_1 \leftarrow \mathbf{L}$
\FOR{$j=1$ to $J$} 
\STATE  $\mathbf{S}_j \leftarrow \underset{\mathbf{S}\in \mathcal{R}_G}{\text{argmin}} \left\Vert \mathbf{S}^T\mathbf{L}_j\mathbf{S} \right\Vert_{\text{offdiag}}^2$
\STATE  $\mathbf{L}_{j+1} \leftarrow \mathbf{S}_j^T\mathbf{L}_j\mathbf{S}_j$
\ENDFOR
\STATE  $\hat{\boldsymbol{\Lambda}} \leftarrow \text{diag}(\mathbf{L}_{J+1})$
\modfinal{\STATE Sort diagonal entries of $\hat{\boldsymbol{\Lambda}}$ in increasing order. Reorder columns of $\mathbf{S}_J$ accordingly.}
\ENSURE sparse orthogonal factors $\mathbf{S}_1,\ldots,\mathbf{S}_J$; diagonal factor $\hat{\boldsymbol{\Lambda}}$.
\end{boxedalgorithmic}
\caption{{\bf Truncated Jacobi algorithm}: Approximate diagonalization algorithm with prescribed complexity.}
\label{algo:givens_sum}
\end{figure}

\subsection{Subproblem resolution}

The algorithm requires to solve $J$ times the optimization subproblem~\eqref{eq:subprob} (at line $3$ of the algorithm of \modLLM{Figure}~\ref{algo:givens_sum}). The solution of this subproblem is given by the Givens rotation $\mathbf{G}_{p,q,\theta}$, where the indices $p$ and $q$ correspond to the greatest entry of $\mathbf{L}_j$ in absolute value (denoted $|l_{pq}^j|$), and the rotation angle has the expression $\theta = \frac{1}{2}\arctan(\frac{l^j_{qq} - l^j_{pp}}{2l^j_{pq}}) + (2k+1)\frac{\pi}{4}$, $k \in \mathbb{Z}$. We then have $\left\Vert \mathbf{L}_{j+1} \right\Vert_{\text{offdiag}}^2 = \left\Vert \mathbf{L}_j \right\Vert_{\text{offdiag}}^2 - 2(l^j_{pq})^2 $. For a proof as well as a review of the different implementations and variants of the Jacobi algorithm, see \cite[pages 426-435]{Golub2012}.  

\begin{figure}[tbp]
\raggedright 
\begin{boxedalgorithmic}[1] 
\REQUIRE{matrix $\mathbf{L}_j$.}
\STATE $(p,q) \leftarrow \underset{(r,s)\in [n]^2}{\text{argmax}} |l^j_{rs}|$ 
\STATE $\theta \leftarrow \frac{1}{2}\arctan(\frac{l^j_{qq} - l^j_{pp}}{2l^j_{pq}}) + \frac{\pi}{4}$
\STATE  $\mathbf{S}_j \leftarrow \mathbf{G}_{p,q,\theta} $
\ENSURE matrix $\mathbf{S}_j = \underset{\mathbf{S}\in \mathcal{R}_G}{\text{argmin}} \left\Vert \mathbf{S}^T\mathbf{L}_j\mathbf{S} \right\Vert_{\text{offdiag}}^2$.
\end{boxedalgorithmic}
\caption{{\bf Resolution of subproblem~\eqref{eq:subprob}}}
\label{algo:subprob}
\end{figure}

\subsection{Parallel truncated Jacobi algorithm} 

 The Relative Complexity Gain RCG, as defined in \modLLM{Section~}\ref{sec:prob_form}, is only a theoretical measure of complexity gain, it corresponds to a ratio between the number of scalar multiplication required to compute the product with the true Fourier matrix $\mathbf{U}$ and the number of scalar multiplications required to compute the product with its efficient approximation $\hat{\mathbf{U}} = \mathbf{S}_1\dots \mathbf{S}_J$. However the actual time gain, that is observed when using the efficient approximation to multiply vectors, is related but not directly proportional to the number of required scalar multiplications. It depends indeed on other factors, such as the implementation of the matrix/vector product, involving or not parallelization. For example, in MATLAB, the product of a vector $\mathbf{x}$ by the dense Fourier matrix $\mathbf{U}$ can be faster than the product with the efficient approximation $\hat{\mathbf{U}} = \mathbf{S}_1\dots \mathbf{S}_J$ even if the RCG is large. This is because the dense matrix/vector product in MATLAB involves parallelism whereas the product with the efficient approximation is done sequentially ($\mathbf{y}\leftarrow\mathbf{S}_J\mathbf{x}$, $\mathbf{y}\leftarrow\mathbf{S}_{J-1}\mathbf{y}$,$\dots$, $\mathbf{y}\leftarrow\mathbf{S}_1\mathbf{y}$). Since some of the Givens rotations could be applied at the same time in parallel (as soon as their supports are disjoint), this leaves room for improvement. 

In order to overcome this limitation, we propose to slightly modify the algorithm of \modLLM{Figure}~\ref{algo:givens_sum} to choose at each step of the approximate diagonalization not only one Givens rotation, but $\frac{n}{2}$ Givens rotations that can be applied in parallel\footnote{For clarity's sake, we suppose here that $n$ is even. If $n$ is odd, 
then $\mathbf{S}$ is composed of $\lfloor\frac{n}{2}\rfloor$ Givens rotation and keeps one coordinate untouched.}. More formally, for an approximation with $J$ Givens rotations, this corresponds to choosing $K \triangleq \lceil\frac{2J}{n}\rceil$ factors, where each of the sparse factors $\mathbf{S}_k$ belongs to the set
$\mathcal{R}_P$ of matrices made of $\frac{n}{2}$ Givens rotations with mutually disjoint supports. 
Elements of this set are matrices of the form
$$
\mathbf{S} = \mathbf{P}
\left(
\begin{array}{ccc}
\mathbf{R}_1&&0 \\
&\ddots&\\
0&&\mathbf{R}_{\frac{n}{2}}\\
\end{array}
\right)
\mathbf{P}^T,
$$
where $\mathbf{P}$ is a permutation matrix and $\forall i \in \{1,\dots,\frac{n}{2}\}$ we have
$$
\mathbf{R}_i = \left(
\begin{array}{cc}
\cos{\theta_i}&-\sin{\theta_i}\\
\sin{\theta_i}&\cos{\theta_i}\\
\end{array}
\right).
$$
The algorithm corresponding to this modified version of the approximate diagonalization is given in \modLLM{Figure}~\ref{algo:givens_parall}. In practice, the choice for the supports of the parallel Givens rotations at the $k$th step is done in a similar way as in the algorithm of \modLLM{Figure}~\ref{algo:givens_sum} (by sequentially choosing the greatest entry of $\mathbf{L}_k$), except that it is necessary to make sure that the support of each Givens rotation be disjoint with the supports of all the Givens rotations chosen previously during the same step. This can be done by:
\begin{enumerate}
\item  Sorting all the nonzero entries of $\mathbf{L}_k$ only once at each step $k$ and put their ordered indices in a list.
\item Choosing $\frac{n}{2}$ Givens rotations to put in the $k$th factor $\mathbf{S}_k$ by going down the list, making sure each chosen Givens rotation has its support disjoint from all the previously chosen ones in the current factor.
\end{enumerate} 
\modfinal{Despite the fact that this strategy is not optimal at each step, $$\mathbf{S}_k \neq \underset{\mathbf{S} \in \mathcal{R}_p}{\text{argmin}} \left\Vert \mathbf{S}^T\mathbf{L}_k\mathbf{S} \right\Vert^2_{\text{offdiag}},$$ it is guaranteed to yield $\left\Vert\mathbf{L}_{k+1} \right\Vert^2_{\text{offdiag}} < \left\Vert\mathbf{L}_{k} \right\Vert^2_{\text{offdiag}}$. As is shown empirically in the experimental part of the paper, this non-optimality does not harm the performance of the method.}

\modfinal{Note that we could also perfectly choose less than $\frac{n}{2}$ rotations to be performed in parallel at each step. This would be interesting if architectural constraints limit the number of arithmetic operations that can be done in parallel ($\frac{n}{2}$ being the optimal choice in case of an architecture with many processors).}

A similar approach where parallel elementary rotations are chosen is evoked in \cite{Kondor2014}, with supplementary constraints due to multiresolution. \textcolor{black}{The approach proposed here relaxes these constraints. This approach is also} different from the parallel MMF (pMMF) factorization method \cite{Teneva2016}, where the main goal is to accelerate the factorization \textcolor{black}{algorithm} itself, by clustering the columns/rows of $\mathbf{L}_k$ every few iterations to reduce the cost of finding the support of the Givens rotations.


\begin{figure}[tbp]
\centering 
\begin{boxedalgorithmic}[1] 
\REQUIRE{matrix $\mathbf{L}$, number $J$ of Givens rotations.}
\STATE $\mathbf{L}_1 \leftarrow \mathbf{L}$
\STATE $j\leftarrow 0$, $k\leftarrow 1$
\WHILE{$j<J$} 
\STATE Choose $\mathbf{S}_k \in \mathcal{R}_P$ {\footnotesize(such that $\left\Vert \mathbf{S}_k^T\mathbf{L}_k\mathbf{S}_k \right\Vert_{\text{offdiag}}^2 < \left\Vert \mathbf{L}_k\right\Vert_{\text{offdiag}}^2$)}
\STATE  $\mathbf{L}_{k+1} \leftarrow \mathbf{S}_k^T\mathbf{L}_k\mathbf{S}_k$
\STATE $j\leftarrow j+\frac{n}{2}$, $k\leftarrow k+1$
\ENDWHILE
\STATE  $\hat{\boldsymbol{\Lambda}} \leftarrow \text{diag}(\mathbf{L}_{k})$, $K \leftarrow k$
\modfinal{\STATE Sort diagonal entries of $\hat{\boldsymbol{\Lambda}}$ in increasing order. Reorder columns of $\mathbf{S}_K$ accordingly.}
\ENSURE sparse orthogonal factors $\mathbf{S}_1,\ldots,\mathbf{S}_K$; diagonal factor $\hat{\boldsymbol{\Lambda}}$.
\end{boxedalgorithmic}
\caption{{\bf Parallel truncated Jacobi algorithm}: Approximate diagonalization algorithm with prescribed complexity and parallel Givens rotations.}
\label{algo:givens_parall}
\end{figure}

\subsection{Computational cost of the approaches}
The proposed approaches yield approximate Fourier transforms with potentially better complexity than the exact Fourier transform obtained by exact diagonalization. The cost of performing the associated (approximate) diagonalization is also different, we detail it below for each approach.

\noindent
 {\bf Exact diagonalization \cite{Jacobi1846}: } exact diagonalization of a symmetric matrix costs $\mathcal{O}(n^3)$ (for example using the classical Jacobi algorithm). \modLLM{In theory, faster diagonalization is possible \cite{Strassen1969,Coppersmith1990,Legall2014}, at the price of higher constants and/or stability issues.}
 Once this diagonalization is computed, and since the obtained eigenvector matrix is in general dense, applying it costs $\mathcal{O}(n^2)$. It is known that the Jacobi algorithm converges at least linearly to a diagonal matrix \cite{Golub2012}.

\noindent {\bf Truncated Jacobi: } the cost of the truncated Jacobi algorithm of \modLLM{Figure}~\ref{algo:givens_sum} is dominated by the quest for the greatest entry of $\mathbf{L}_j$ (line $1$ of the algorithm of \modLLM{Figure}~\ref{algo:subprob}), that costs $\mathcal{O}(\|\mathbf{L}_{j}\|_{0})$ (which is at worst $\mathcal{O}(n^2)$, but can be much smaller if $\mathbf{L}_{j}$ is sparse). Since matrices $\mathbf{L}_j$ and $\mathbf{L}_{j+1}$ have identical entries except for the lines and columns $p$ and $q$, one can reuse computations carried out at the previous step for the current step, using a technique inspired by \cite{Cao2011} (a full algorithm implementing this technique is given in appendix~\ref{app:efficientdiago}). Using this technique, the worst case complexity of the subproblem resolution remains \textcolor{black}{$\mathcal{O}(\|\mathbf{L}_{j}\|_{0})$}, 
but drops down to $\mathcal{O}(n)$ for most iterations (as soon as the selected coordinates $p$ and $q$ at the current iteration are both different from the two selected at the previous iteration). This results in a complete truncated Jacobi costing $\mathcal{O}(\textcolor{black}{n^{2}+}nJ)$ arithmetic operations in average. 

Once the approximate diagonalization is done, applying the resulting $\hat{\mathbf{U}}$ costs $\mathcal{O}(J)$ arithmetic operations, since applying a Givens rotation costs only $4$ multiplications and $2$ additions.

\noindent {\bf Parallel truncated Jacobi: } 
its cost is dominated by the sorting done at each step $k$. Since there are at most $\tfrac{n(n-1)}{2}$ entries to sort at each step, the cost of the sorting algorithm can be assumed to be $\mathcal{O}(n^2 \log n)$ 
(the average complexity of most sorting algorithms with $r$ elements being $\mathcal{O}(r\log r)$). The sorting has to be done $K = \mathcal{O}(\frac{J}{n})$ times, which brings the overall complexity of the greedy diagonalization strategy of \modLLM{Figure}~\ref{algo:givens_parall} to $\mathcal{O}(nJ\log n)$. 

As with the plain truncated Jacobi case, once the approximate diagonalization is done, applying the resulting $\hat{\mathbf{U}}$ costs $\mathcal{O}(J)$ arithmetic operations. Yet, since $\frac{n}{2}$ Givens rotations have disjoint supports in each sparse factor, computations can be parallelized leading to much faster approximate \modLLM{FGFTs} in practice \textcolor{black}{as we will see in \modLLM{Section~}\ref{sec:sensor}}. 

\arraycolsep=1pt\def\arraystretch{1.25}
\begin{table}[tb]                                                       
\centering                                                              
\begin{tabularx}{1\columnwidth}{Xcc}                                      
\toprule                                                       
&  Obtaining $\mathbf{U}$ or $\hat{\mathbf{U}}$   & Applying $\mathbf{U}$ or $\hat{\mathbf{U}}$    \\ 
\midrule 
Exact diagonalization \cite{Jacobi1846} & $\mathcal{O}(n^3)$ &$\mathcal{O}(n^2)$  \\
\hdashline[0.5pt/2pt]
Truncated Jacobi & $\mathcal{O}(n^2 +nJ)$ &$\mathcal{O}(J)$ \\
\---- with $J=\mathcal{O}(n \log n)$ & $\mathcal{O}(n^2\log n)$ &$\mathcal{O}(n \log n)$ \\
 \hdashline[0.5pt/2pt]
Parallel truncated Jacobi & $\mathcal{O}(nJ\log n)$ &$\mathcal{O}(J)$ \\
\---- with $J=\mathcal{O}(n \log n)$ & $\mathcal{O}(n^2\log^2 n)$ &$\mathcal{O}(n \log n)$ \\
\bottomrule                                                             
\end{tabularx}                                                           
~\vspace{0mm}\caption{Complexity of exact {\em vs} approximate diagonalizations.
}                                              
\label{tab:complexity}                                              
\end{table} 

\noindent {\bf Summary:}
A comparison between the computational complexities of the classical method to obtain the graph Fourier transform (an exact diagonalization of the Laplacian) and the approximate diagonalization methods proposed here is given in \modLLM{Table~}\ref{tab:complexity}. For the approximate diagonalizations we upper bound $\|\mathbf{L}_{j}\|_{0}$ and  $\|\mathbf{L}_{k}\|_{0}$ by $\mathcal{O}(n^2)$, since even if in the beginning the Laplacian is sparse, applying Givens rotations rapidly introduces nonzero entries. This comparison indicates that in order for the approximate diagonalization to be cheaper than the exact diagonalization, one should take $J = \mathcal{O}(n^\alpha)$ with $\alpha <2$, or even $J = \mathcal{O}(n\log n)$ (which leads to a cost of application of the same order as the classical FFT).


Interestingly, for general symmetric matrices, when diagonalization is done with the Jacobi algorithm, it has been empirically observed that $\mathcal{O}(n^2\log n)$ Givens rotations are typically required to achieve a prescribed tolerance \cite{Golub2012,Brent1985}. One of the main contributions of this paper is to show empirically in the next section that, for certain graph Laplacians, the truncated algorithms (of \modLLM{Figure}s~\ref{algo:givens_sum} and \ref{algo:givens_parall}) provide a good approximate diagonalization with only $J = \mathcal{O}(n\log n)$ Givens rotations, thus enabling an efficient approximate \modLLM{FGFT} on the considered graph.


\section{Evaluation of the approximate Fourier transforms}
\label{ssec:eval_random}
\label{sec:expval}

In this section, we compare the approximate graph Fourier transforms obtained with the methods introduced in this paper to the factorization methods of \cite{Lemagoarou2016}. We first introduce the various graphs we use for the evaluation. We then present the different approximate \modLLM{FGFT} methods we will compare, as well as the different performance measures. We finish by discussing the obtained results. \modfinal{A toolbox implementing all the algorithms and experiments performed in this paper is available at \texttt{http://faust.inria.fr}, in the folder demos\_and\_applications/FGFT.}

\modLLM{\subsection{General evaluation}}

\begin{figure}[tb]
\centering
\includegraphics[width=\columnwidth]{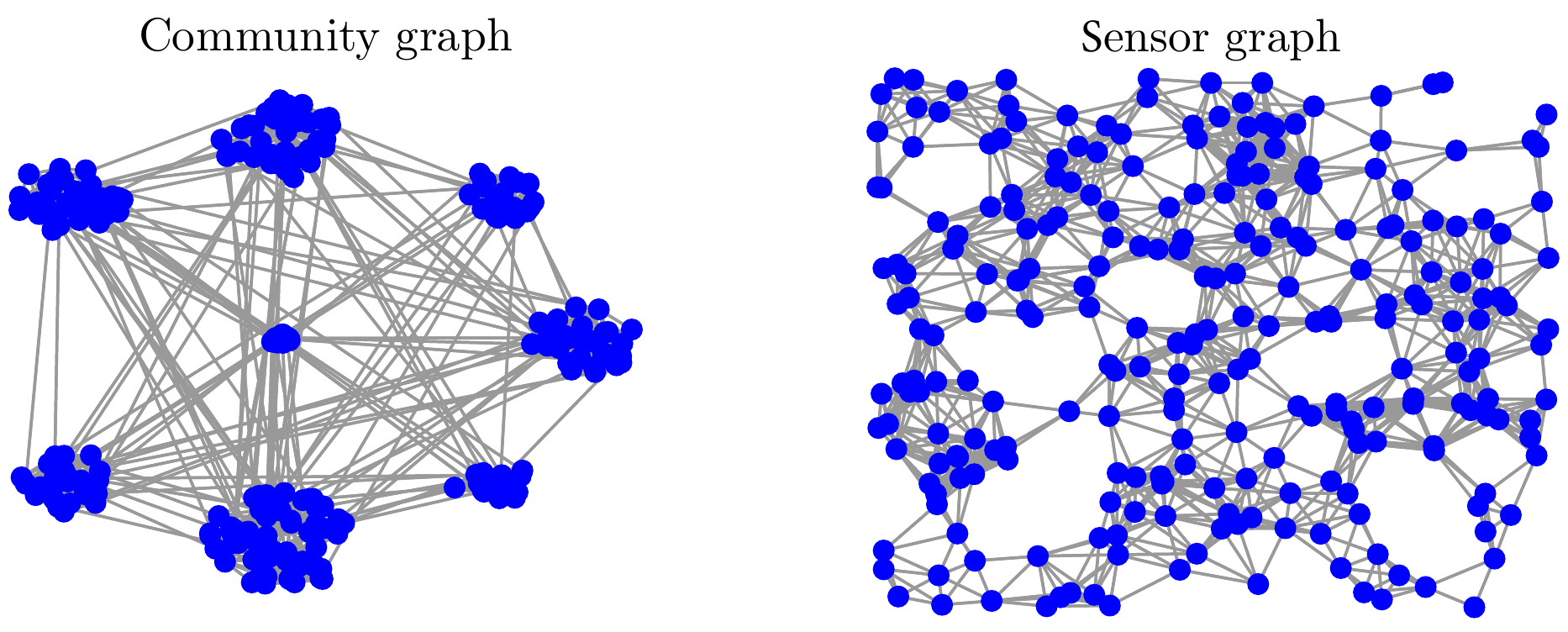}
\caption{Examples of  a community graph and a sensor graph, of size $n=256$.} 
\label{fig:graphexample}
\end{figure}

\noindent {\bf Considered graphs.} For this experiment, we consider several families of graphs among the most common ones. All graphs used here are generated with the help of the ``Graph Signal Processing'' toolbox (GSPBOX) \cite{GSP}. 
\begin{itemize}
\item {\bf Erd\H{o}s-R\'enyi:} a totally random graph where every pair of nodes is connected with a unit weight with probability $p=0.1$.
\item {\bf Community:} a graph of size $n$ made of $\sqrt{n}/2$ communities of random sizes. Every community is itself a graph whose nodes correspond to points whose location is drawn uniformly at random on the unit disk, that are connected with a weight inversely proportional to their distance, if this distance is inferior to a certain threshold. Moreover, random weights between any pair of nodes appear with a probability of $1/n$.
\item {\bf Sensor:} a graph of size $n$ where the nodes correspond to points whose location is drawn uniformly at random on the unit square, that are connected with a weight inversely proportional to their distance, if this distance is inferior to a certain threshold, zero otherwise. 
\modLLM{\item {\bf Ring:} the classical ring graph, where each node is connected to the two adjacent nodes with edges of equal weights.}
\end{itemize}
For all these families of graphs, we take graphs of various sizes $n$, with $n\in \{128,256,512,1024\}$ nodes. \modLLM{As our purpose here is to qualitatively compare the approximate diagonalization method with the factorization methods of \cite{Lemagoarou2016} (which do not scale well to large graphs, as they involve a full diagonalization), we had to restrict to these graph sizes that may seem small. We expect the conclusions of this section to hold for much bigger graphs}. Examples of graphs used in these experiments are shown on \modLLM{Figure}~\ref{fig:graphexample}.

\noindent {\bf Computation of the Fourier transforms.} We consider for each configuration four approximate graph Fourier transforms of complexity $\mathcal{O}(n \log n)$, thus mimicking the complexity of the classical FFT. The corresponding relative complexity gain RCG is of order $\mathcal{O}(n/\log n)$, and is the same for the all the approximate transforms. 
\begin{itemize}
 \item {\bf $\hat{\mathbf{U}}_{\text{PALM}}$:} computed with the factorization method of \cite{Lemagoarou2016}, where the true Fourier matrix is hierarchically factorized into $K$ factors (sparse but not orthogonal; sparsity levels described \modLLM{in detail} in \cite{Lemagoarou2016}) using the Proximal Alternating Linearized Minimization (PALM) algorithm \cite{Bolte2014}.
 \item {\bf $\hat{\mathbf{U}}_{\text{PALM}_d}$:} computed with the version of the factorization method of \cite{Lemagoarou2016} whose global optimization steps target a diagonalization cost function (aiming at an efficient approximation of the Fourier transform that approximately diagonalizes the Laplacian matrix). 
 \item {\bf $\hat{\mathbf{U}}_{\text{Givens}}$:} computed with the truncated Jacobi algorithm. 
 \item {\bf $\hat{\mathbf{U}}_{\text{Givens}\scriptscriptstyle //}$:} 
 computed with the parallel truncated Jacobi algorithm.
\end{itemize}

\noindent {\bf Performance measures.} How can we measure the \modLLM{accuracy} of an approximate graph Fourier transform $\hat{\mathbf{U}}$? 
For this experiment we choose two performance measures. The first one is the relative error on the calculation of the Fourier basis:
$$
\text{err}_c(\hat{\mathbf{U}}) \triangleq \frac{\Vert \mathbf{U} - \hat{\mathbf{U}} \Vert_F}{\left\Vert \mathbf{U} \right\Vert_F}.
$$
\modfinal{Note that for the approximate FGFTs obtained via diagonalization of the Laplacian ($\hat{\mathbf{U}}_{\text{Givens}}$ and $\hat{\mathbf{U}}_{\text{Givens}\scriptscriptstyle //}$), the sign of each of their columns is adjusted so that its scalar product with the corresponding column of $\mathbf{U}$ is positive. This procedure allows to get rid of the sign ambiguity inherent to the diagonalization methods before the computation of $\text{err}_c(\hat{\mathbf{U}})$.}

The second performance measure is the relative error on the diagonalization of the Laplacian:
$$
\text{err}_d(\hat{\mathbf{U}}) \triangleq \frac{\Vert \hat{\mathbf{U}}^T \mathbf{L} \hat{\mathbf{U}}\Vert_{\text{offdiag}}}{\left\Vert \mathbf{L} \right\Vert_F}.
$$

\arraycolsep=1pt\def\arraystretch{1.25}
\begin{table*}[tb]                                                       
\centering                                                              
\begin{tabularx}{1\textwidth}{Xccccc}                                      
\toprule                                                       
& & Erd\H{o}s-R\'enyi  & Community & Sensor & Ring \\
& &
\resizebox{27mm}{!}{$\hat{\mathbf{U}}_\text{PALM} | \hat{\mathbf{U}}_{\text{PALM}_d} | \hat{\mathbf{U}}_\text{Givens} | \hat{\mathbf{U}}_{\text{Givens}\scriptscriptstyle //}$}
&
\resizebox{27mm}{!}{$\hat{\mathbf{U}}_\text{PALM} | \hat{\mathbf{U}}_{\text{PALM}_d} | \hat{\mathbf{U}}_\text{Givens} | \hat{\mathbf{U}}_{\text{Givens}\scriptscriptstyle //}$}
&
\resizebox{27mm}{!}{$\hat{\mathbf{U}}_\text{PALM} | \hat{\mathbf{U}}_{\text{PALM}_d} | \hat{\mathbf{U}}_\text{Givens} | \hat{\mathbf{U}}_{\text{Givens}\scriptscriptstyle //}$}
&
\resizebox{27mm}{!}{$\hat{\mathbf{U}}_\text{PALM} | \hat{\mathbf{U}}_{\text{PALM}_d} | \hat{\mathbf{U}}_\text{Givens} | \hat{\mathbf{U}}_{\text{Givens}\scriptscriptstyle //}$}\\
\midrule                                                                 
\multirow{3}{*}{
\small $n=128$, $\text{RCG} = 3.9$ 
}  
 &$\text{err}_c(\hat{\mathbf{U}})$ &{\bf 0.44} \textbar 0.53 \textbar 1.23 \textbar 1.25 & {\bf 0.06} \textbar 0.11 \textbar 0.90 \textbar 1.01 & {\bf 0.17} \textbar 0.22 \textbar 0.85 \textbar 0.97 &{\bf 0.53} \textbar 0.65 \textbar 1.18 \textbar 1.17 \\
&$\text{err}_d(\hat{\mathbf{U}})$ &0.45 \textbar 0.38 \textbar {\bf 0.11} \textbar 0.12  &0.08 \textbar 0.07 \textbar {\bf 0.03} \textbar 0.05 & 0.21 \textbar 0.19 \textbar {\bf 0.05} \textbar 0.08 & 0.39 \textbar 0.45 \textbar {\bf 0.09} \textbar {\bf 0.09} \\
\cdashline{2-6}[.5pt/2pt]  
& $T$ & 6.63 \textbar 5.22 \textbar 0.40 \textbar {\bf 0.14} & 7.24 \textbar 5.44 \textbar 0.44 \textbar {\bf 0.13} & 7.00 \textbar 5.33 \textbar 0.44 \textbar {\bf 0.11} & 7.03 \textbar 5.44 \textbar 0.39 \textbar {\bf 0.10} \\
\midrule                                                                  
\multirow{3}{*}{
\small $n=256$, $\text{RCG} = 7.1$
}
&$\text{err}_c(\hat{\mathbf{U}})$ & {\bf 0.61} \textbar 0.86 \textbar 1.31 \textbar 1.31 & {\bf 0.19} \textbar 0.33 \textbar 1.08 \textbar 1.14 & {\bf 0.27} \textbar 0.36 \textbar 1.07 \textbar 1.15 &{\bf 0.69} \textbar 0.79 \textbar 1.24 \textbar 1.24 \\ 
&$\text{err}_d(\hat{\mathbf{U}})$& 0.45 \textbar 0.34 \textbar {\bf 0.10} \textbar 0.11 & 0.21 \textbar 0.18 \textbar {\bf 0.05} \textbar 0.07 & 0.30 \textbar 0.24 \textbar {\bf 0.06} \textbar 0.08 & 0.40 \textbar 0.50 \textbar {\bf 0.08} \textbar {\bf 0.08} \\
\cdashline{2-6}[.5pt/2pt]  
&$T$ & 56.4 \textbar 47.9 \textbar 1.60 \textbar {\bf 0.60} & 60.8 \textbar 47.1 \textbar 1.52 \textbar {\bf 0.52} & 59.0 \textbar 45.2 \textbar 1.67 \textbar {\bf 0.41} & 55.8 \textbar 44.7 \textbar 1.59 \textbar {\bf 0.27} \\
\midrule  
\multirow{3}{*}{
\small $n=512$, $\text{RCG} = 13.1$
}
&$\text{err}_c(\hat{\mathbf{U}})$& {\bf 0.73} \textbar 1.00 \textbar 1.35 \textbar 1.35 & {\bf 0.30} \textbar 0.85 \textbar 1.20 \textbar 1.25 & {\bf 0.34} \textbar 0.50 \textbar 1.20 \textbar 1.25 & {\bf 0.81} \textbar 0.94 \textbar 1.30 \textbar 1.29 \\ 
&$\text{err}_d(\hat{\mathbf{U}})$&0.38 \textbar {\bf 0.05} \textbar 0.08 \textbar 0.08 & 0.27 \textbar 0.17\textbar {\bf 0.05} \textbar 0.07 & 0.33 \textbar 0.28 \textbar {\bf 0.06} \textbar 0.08 &0.34 \textbar 0.61 \textbar {\bf 0.07} \textbar 0.08 \\
\cdashline{2-6}[.5pt/2pt] 
&$T$ & 498 \textbar 461 \textbar 14.9 \textbar {\bf 5.18} & 511 \textbar 459 \textbar 14.6 \textbar {\bf 4.37} & 515 \textbar 457 \textbar 15.0 \textbar {\bf 3.06} & 494 \textbar 453 \textbar 14.4 \textbar {\bf 1.33}  \\
\midrule  
\multirow{3}{*}{
\small $n=1024$, $\text{RCG} = 24.4$
}
&$\text{err}_c(\hat{\mathbf{U}})$& {\bf 0.82} \textbar 1.00 \textbar 1.37 \textbar 1.37 & {\bf 0.42} \textbar 1.00 \textbar 1.28 \textbar 1.31 & {\bf 0.39} \textbar 0.64 \textbar 1.29 \textbar 1.32 & {\bf 0.87} \textbar 1.07 \textbar 1.33 \textbar 1.32 \\ 
&$\text{err}_d(\hat{\mathbf{U}})$& 0.29 \textbar {\bf 0.02} \textbar 0.06 \textbar 0.07 & 0.31 \textbar 0.07 \textbar {\bf 0.05} \textbar 0.06 & 0.35 \textbar 0.31 \textbar {\bf 0.06} \textbar 0.08 & 0.28 \textbar 0.77 \textbar {\bf 0.07} \textbar {\bf 0.07} \\ 
\cdashline{2-6}[.5pt/2pt] 
&$T$ & 6e+3 \textbar 5e+3 \textbar 122 \textbar {\bf 52.2} & 6e+3 \textbar 5e+3 \textbar 119 \textbar {\bf 42.0} & 6e+3 \textbar 5e+3 \textbar 120 \textbar {\bf 24.9} & 6e+3 \textbar 5e+3 \textbar 117 \textbar {\bf 7.77} \\
\bottomrule                                                             
\end{tabularx}                                                           
~\vspace{0mm}\caption{Results of the evaluation of the approximate fast Fourier transforms. 
For each configuration, the results are given for the four methods evoked in the text
, and the best result is written in bold. 
}                                                
\label{tab:resultgraphs}                                              
\end{table*}

\noindent {\bf Results.} Average results (over $10$ realizations) are shown in \modLLM{Table~}\ref{tab:resultgraphs}. 
In addition to the performance measures, we also show the factorization time $T$ (in seconds, not taking into account the \textcolor{black}{$\mathcal{O}(n^{3})$} time needed to obtain $\mathbf{U}$ for the first two methods). Several comments are in order: 
\begin{itemize}
\item First, it is very clear that the factorization method of \cite{Lemagoarou2016} gives approximate \modLLM{FGFT}s $\hat{\mathbf{U}}_\text{PALM}$ that are better in terms of relative calculation error $\text{err}_c$, compared to all the other methods tested here. Conversely, the diagonalization methods presented here gives approximate \modLLM{FGFT}s $\hat{\mathbf{U}}_\text{Givens}$ and $\hat{\mathbf{U}}_{\text{Givens}\scriptscriptstyle //}$ that are better in terms of relative diagonalization error $\text{err}_d$. This is true for almost all the tested configurations (except for the Erd\H{o}s-R\'enyi graph in dimension $512$ and $1024$ where the approximate \modLLM{FGFT} $\hat{\mathbf{U}}_{\text{PALM}_d}$ gives better diagonalization results). This seems quite logical and is concordant with the cost functions respectively considered by the methods. \modLLM{Indeed, the cost function used for the factorization methods is very close to $\text{err}_c$ and the cost function used for the diagonalization methods is very close to $\text{err}_d$.} Moreover, performance of $\hat{\mathbf{U}}_\text{Givens}$ and $\hat{\mathbf{U}}_{\text{Givens}\scriptscriptstyle //}$ are very close, showing that taking parallel Givens rotations does not decrease performance, while providing approximate \modLLM{FGFT}s that can be \modLLM{applied (and obtained)} much faster. 
\item Second, regarding the different graph families considered in the experiment, we see that all methods show in general better results for the graphs ``Sensor'' and ``Community'', and poorer results for the graphs ``Erd\H{o}s-R\'enyi'' and ``Ring''. It is quite expected that it is difficult to obtain a good approximate \modLLM{FGFT} for the graph ``Erd\H{o}s-R\'enyi'', since it is totally random and unstructured. However, it is a bit more surprising for the graph ``Ring'', that is highly structured, and for which we know a fast Fourier transform exists \modLLM{(which is nothing more than the classical butterfly FFT \cite{CooleyTukey1965}).}
This indicates that the optimization techniques used to obtain the approximate \modLLM{FGFT}s (both the factorizations and the greedy diagonalizations) are in a way suboptimal, in the sense that they do not attain a global optimum of their cost functions. \modLLM{We discuss links between the graph structure and the diagonalization performance in section~\ref{ssec:structure_perf}.}
Moreover, the error measure \modLLM{$\text{err}_c$} is greater than one in certain configurations for the approximate \modLLM{FGFT}s $\hat{\mathbf{U}}_{\text{PALM}_d}$, $\hat{\mathbf{U}}_\text{Givens}$ and $\hat{\mathbf{U}}_{\text{Givens}\scriptscriptstyle //}$. This can be explained by the fact that the cost \modLLM{functions} for these methods do not enforce the true Fourier matrix $\mathbf{U}$ and its approximation to be close to each other in Frobenius norm.
\item Third, we notice that for $\hat{\mathbf{U}}_\text{Givens}$ and $\hat{\mathbf{U}}_{\text{Givens}\scriptscriptstyle //}$, $\text{err}_d$ does not increase with the graph size $n$, although the relative complexity gain increases (we consider here approximate \modLLM{FGFT}s whose complexity is $\mathcal{O}(n\log n)$). In other words, the computational benefit grows with $n$ when using the approximate \modLLM{FGFT}s, without increasing the relative error: this behavior is typically what is expected from an approximate fast transform. The same behavior was observed~\cite{Lemagoarou2016} for the factorization method, with approximate \modLLM{FGFT}s of higher complexity $\mathcal{O}(n^{1.27})$.
\item Last, regarding the factorization time $T$, $\hat{\mathbf{U}}_\text{Givens}$ and $\hat{\mathbf{U}}_{\text{Givens}\scriptscriptstyle //}$ are obtained much faster (often by two orders of magnitude) than  $\hat{\mathbf{U}}_\text{PALM}$ and $\hat{\mathbf{U}}_{\text{PALM}_d}$. The truncated Jacobi algorithms presented in this paper are thus faster than the methods proposed in  \cite{Lemagoarou2016}. Moreover the parallel truncated Jacobi algorithm of \modLLM{Figure}~\ref{algo:givens_parall} (that gives $\hat{\mathbf{U}}_{\text{Givens}\scriptscriptstyle //}$) is faster than the plain truncated Jacobi algorithm of \modLLM{Figure}~\ref{algo:givens_sum} (that gives $\hat{\mathbf{U}}_\text{Givens}$), \textcolor{black}{although its theoretical complexity is higher (see \modLLM{Table~}\ref{tab:complexity}). This may be due to the respective implementations of the two methods, or to the constants before the complexity orders being different}. For example, it is three to five times faster depending on the graph type for $n=1024$.
\end{itemize}

In summary, the approximate \modLLM{FGFT}s $\hat{\mathbf{U}}_\text{PALM}$ and $\hat{\mathbf{U}}_{\text{PALM}_d}$ computed by the factorization method are more adapted to situations where the knowledge of the actual Fourier coefficients of a signal $\mathbf{x}$ is important ($\mathbf{U}^T\mathbf{x}$ and $\hat{\mathbf{U}}^T\mathbf{x}$ are close). Conversely, the approximate \modLLM{FGFT}s $\hat{\mathbf{U}}_{\text{Givens}\scriptscriptstyle //}$ or $\hat{\mathbf{U}}_{\text{Givens}}$ computed by the diagonalization methods of \modLLM{Figure}s~\ref{algo:givens_sum} or \ref{algo:givens_parall} are more adapted to situations where the ``physical'' interpretation of the Fourier transform of a signal $\mathbf{x}$ is important ($\hat{\mathbf{U}}$ is an orthogonal matrix whose columns are approximately eigenvectors of the Laplacian $\mathbf{L}$). However, the algorithm of \modLLM{Figure}~\ref{algo:givens_parall} shows similar performance but is faster than the algorithm of \modLLM{Figure}~\ref{algo:givens_sum}, and it gives approximate \modLLM{FGFT}s that can be applied faster: it is thus to be preferred in applications. 
Moreover, all methods seem to perform better on structured graphs (``Sensor'' and ``Community'') that are anyway more likely to be encountered in real-world applications.

\noindent {\bf Note.} 
\modNT{Hereafter, and given these first results,  we concentrate on the performance analysis of the approximate \modLLM{FGFT} $\hat{\mathbf{U}}_{\text{Givens}\scriptscriptstyle //}$ computed with the algorithm of \modLLM{Figure}~\ref{algo:givens_parall}.} 
%

\modLLM{
\subsection{Link between structure and diagonalization performance}
\label{ssec:structure_perf}}

\begin{table}[tb]                                       
\modLLM{\centering   
\def\arraystretch{1.2}                                           
\begin{tabular}{rcccc} 
\toprule                                                                              
 &$\epsilon=\tfrac{\epsilon_c}{100}$  & $\epsilon=\tfrac{\epsilon_c}{25}$  &$\epsilon=\tfrac{\epsilon_c}{10}$&$\epsilon=\tfrac{\epsilon_c}{2}$\\
 & \resizebox{10mm}{!}{$\text{err}_d $ \textbar $\text{err}_s $} & \resizebox{10mm}{!}{$\text{err}_d $ \textbar $\text{err}_s $} & \resizebox{10mm}{!}{$\text{err}_d $ \textbar $\text{err}_s $} & \resizebox{10mm}{!}{$\text{err}_d $ \textbar $\text{err}_s $}\\
\midrule                                                                                                           
$\bar{d}=4$&0.036 \textbar 0.002 & 0.061 \textbar 0.004 & 0.090 \textbar 0.008 & 0.140 \textbar 0.019 \\                                                 
$\bar{d}=8$&0.043 \textbar 0.003 & 0.075 \textbar 0.007 & 0.111 \textbar 0.015 & 0.164 \textbar 0.034 \\                                                 
$\bar{d}=16$&0.041 \textbar 0.003 & 0.074 \textbar 0.010 & 0.105 \textbar 0.019 & 0.148 \textbar 0.037 \\                                                 
$\bar{d}=32$&0.037 \textbar 0.004 & 0.065 \textbar 0.012 & 0.089 \textbar 0.021 & 0.120 \textbar 0.034 \\
\bottomrule                                                                                                      
\end{tabular}                                           
~\vspace{1mm}\caption{Relative diagonalization error $\text{err}_d$ \modfinal{and relative spectrum error $\text{err}_s$} on a Stochastic Block Model (SBM) graph, in function of the average degree $\bar{d}$ and the \modNT{structural} parameter $\epsilon$.}                                
\label{tab:structure} }                             
\end{table}

\modLLM{As mentioned in the previous subsection, it seems difficult to establish a direct link between the graph structure and the results obtained by the different methods in terms of diagonalization and factorization error. Focusing on the approximate FGFT $\hat{\mathbf{U}}_{\text{Givens}\scriptscriptstyle //}$ and the performance measure $\text{err}_d$, we perform here approximate diagonalizations on graphs whose level of structure is parameterized in order to gain insight on its influence on the diagonalization performance. \modfinal{Moreover, we take advantage of this subsection to evaluate the approximate eigenvalues computed by the method, via the measure $$\text{err}_s(\hat{\boldsymbol{\Lambda}}) \triangleq \frac{\Vert \text{diag}(\boldsymbol{\Lambda} - \hat{\boldsymbol{\Lambda}})\Vert_2}{\Vert \text{diag}(\boldsymbol{\Lambda})\Vert_2},$$ that is the relative error on the graph spectrum.} To this end, we take Stochastic Block Model (SBM) graphs (see for instance \cite[Section~5.1.]{Tremblay2016}) with $n=1000$ nodes made of twenty communities of equal sizes. We fix the relative complexity gain for this experiment to $\text{RCG}=10$. We vary the average degree $\bar{d}$ of the graph and the ratio $\epsilon$ between in-community and inter-community connection probabilities. Results in average over ten realizations are given in Table~\ref{tab:structure}, as a function of $\bar{d}$ and $\epsilon$ ($\epsilon_c$ is the threshold above which community detection is impossible~\cite{decelle_asymptotic_2011}). First, the average degree does not seem to influence the diagonalization performance. On the other hand, the parameter $\epsilon$ seems to influence the result, since the closer it is to $\epsilon_c$, the worst are the diagonalization results. In other words, more structured SBM graphs (with a low $\epsilon$) lead to better diagonalization performance. In contrast with the previous subsection where it was difficult to interpret the result for different graph families, we see here that in a specific graph family such as SBM graphs, global structure is linked to the diagonalization performance.}

\modfinal{
Finally, regarding the eigenvalues, it can be seen that their approximation is good for every configuration. Indeed, the relative spectrum error $\text{err}_s$ does not exceed $0.037$, compared to the relative diagonalization error which goes up to $0.164$. Among the different tested configurations, eigenvalues are better approximated when the average degree $\bar{d}$ is low (when there are less edges). They are also better approximated with more structured graphs (when $\epsilon$ is low).
}\\

\section{Application to sensor networks}
\label{sec:sensor}
In this section, we first discuss the relevance of our method for graphs representing sensor networks, before analyzing its performance on both synthetic and real-world sensor networks. 

\begin{table*}[tb]                                              
\centering                                                     
\begin{tabular}{rcccccccc}                                   
\toprule                                                       
 & $n=64$ & $n=128$ & $n=256$ & $n=512$ & $n=1024$ & $n=2048$ & $n=4096$ & $n=8192$  \\
\midrule                                                          
RCG & 1.33 & 2.29 & 4.00 & 7.11 & 12.80 & 23.27 & 42.67 & 78.77 \\                        
                                                         
Time gain & 0.04 & 0.05& 0.11 & 0.26 & 1.56 & 3.88 & 7.57 & 27.16 \\                        
$\text{err}_d(\hat{\mathbf{U}}_{\text{Givens}\scriptscriptstyle //})$ & 0.057 &0.062 & 0.055 &0.051 & 0.048 & 0.051 & 0.049 & 0.048 \\                                                           
\bottomrule                                                       
\end{tabular}                                                  
~\vspace{1mm}\caption{Relative Complexity Gain, actual time gain and relative error results on random sensor graphs of various sizes.}
\label{tab:rand_sensor}                                 
\end{table*}

\subsection{Why sensor networks ?}
To obtain the approximate \modLLM{FGFT} $\hat{\mathbf{U}}_{\text{Givens}\scriptscriptstyle //}$, the parallel truncated Jacobi algorithm of \modLLM{Figure}~\ref{algo:givens_parall} costs $\mathcal{O}(nJ \textcolor{black}{\log n})$ arithmetic operations. For instance, it costs $\mathcal{O}(n^2\textcolor{black}{(}\log n\textcolor{black}{)^{2}})$ for $J = \mathcal{O}(n \log n)$ Givens rotations.
This is much 
 cheaper than an exact diagonalization, which costs $\mathcal{O}(n^{3})$, but substantially more costly than applying the resulting approximate fast transform $\mathcal{O}(n \log n)$. 

Favorable use cases for this method are thus applications where this \textcolor{black}{overhead} cost can pay off. This corresponds to cases where the same graph Fourier transform is used a great number of times. Said otherwise, this corresponds to cases where the graph is relatively constant over time, \textcolor{black}{while we process several changing signals over this graph}.

Some graphs are changing relatively fast by nature. For example, graphs corresponding to social network friendships, internet links, or movie ratings are of this kind. These types of applications are thus not well suited for the method, since recomputing an approximate diagonalization each time the graph changes would be too costly.

On the other hand, graphs corresponding to sensors in physical infrastructures, such as road networks or power grids are in general quite constant over time. This kind of graph is thus particularly well suited to the method proposed in this paper, since an approximate diagonalization stays valid for a long time (structural changes in a road network or a power grid do not occur very often). In this context, it can be useful to have knowledge (even if only approximate) of the spectrum of graph signals, for example for monitoring purposes or malfunction detection, see e.g.\ \cite[section V.A]{Sandryhaila2014a}.

\subsection{Experiments on random sensor graphs}

We 
perform approximate diagonalizations of Laplacians corresponding to random sensor graphs  generated with the GSPBOX \cite{GSP}. For this experiment, we take graphs of various sizes $n \in \{64,128,256,512,1024,2048,4096,8192\}$, 
and choose $J=2n\log n$ Givens rotations,  
yielding an approximate Fourier transform of
 the same complexity as the famous classical butterfly FFT \cite{CooleyTukey1965}. Once an approximate \modLLM{FGFT} $\hat{\mathbf{U}}_{\text{Givens}\scriptscriptstyle //}$ is obtained this way, its actual time gain is obtained by measuring the mean time taken in MATLAB to perform matrix/vector products $\hat{\mathbf{U}}_{\text{Givens}\scriptscriptstyle //}\mathbf{x}$ with random vectors $\mathbf{x} \in \mathbb{R}^{n}$, compared to the time taken to perform the matrix/vector products $\mathbf{U}\mathbf{x}$ (the mean being taken over $100$ such products). Results for this experiment are given in \modLLM{Table~}\ref{tab:rand_sensor}, in terms of Relative Complexity Gain, actual time gain and relative diagonalization error $\text{err}_d(\hat{\mathbf{U}}_{\text{Givens}\scriptscriptstyle //})$. 
Several comments are in order:
\begin{itemize}
\item First of all, the larger the graph, the higher the actual computational gains. Indeed, the actual time gain is very poor when $n$ is small ($0.04$ for $n=64$, which means the approximate \modLLM{FGFT} is more than twenty times {\em slower} than the true Fourier transform), become greater than one for $n=1024$, and attains $27.16$ when $n=8192$ (the approximate \modLLM{FGFT} is then more than $27$ times {\em faster} than the true Fourier transform). Moreover, the ratio between RCG and actual time gain decreases as $n$ grows (it is equal to $46$ for $n=128$ and goes down to $3$ for $n=8192$). This means that the RCG is an overly optimistic estimation of the actual time gain for small $n$, but becomes of the right order of magnitude for large $n$.
\item Second, the relative diagonalization error $\text{err}_d$ is approximately constant (around $0.05$) and does not seem to depend on the graph size. This indicates that for random sensor graphs, approximate \modLLM{FGFT}s of $\mathcal{O}(n\log n)$ complexity make perfect sense. The compromise between error and complexity could indeed be controlled by a multiplicative constant $C$, taking $J = C .n \log n$ Givens rotations for the approximate diagonalization. 
\end{itemize}

\begin{figure*}[tb]
  \centering
  \begin{tabular}{ccc}
    \includegraphics[width=0.31\textwidth]{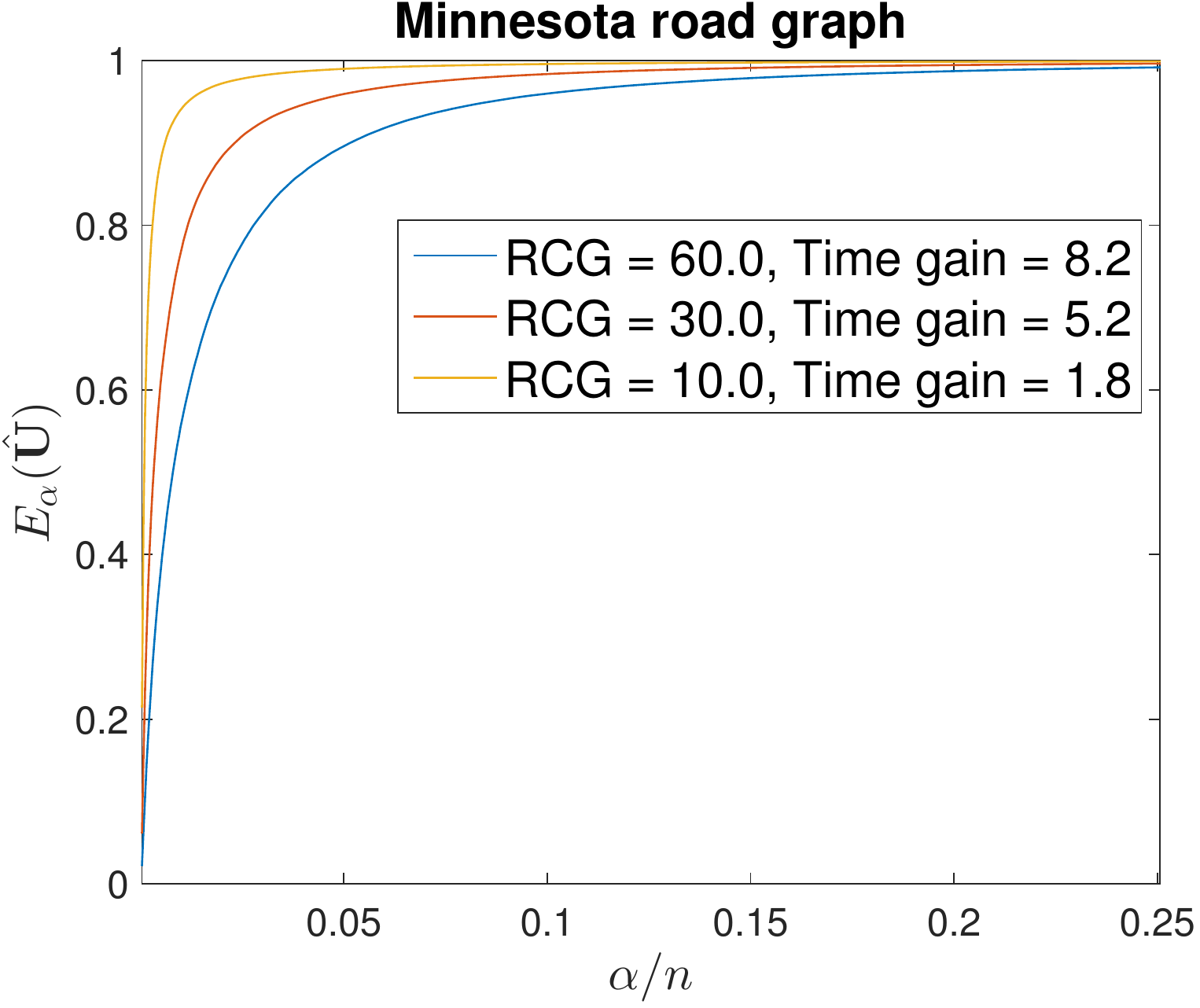}&

    \includegraphics[width=0.31\textwidth]{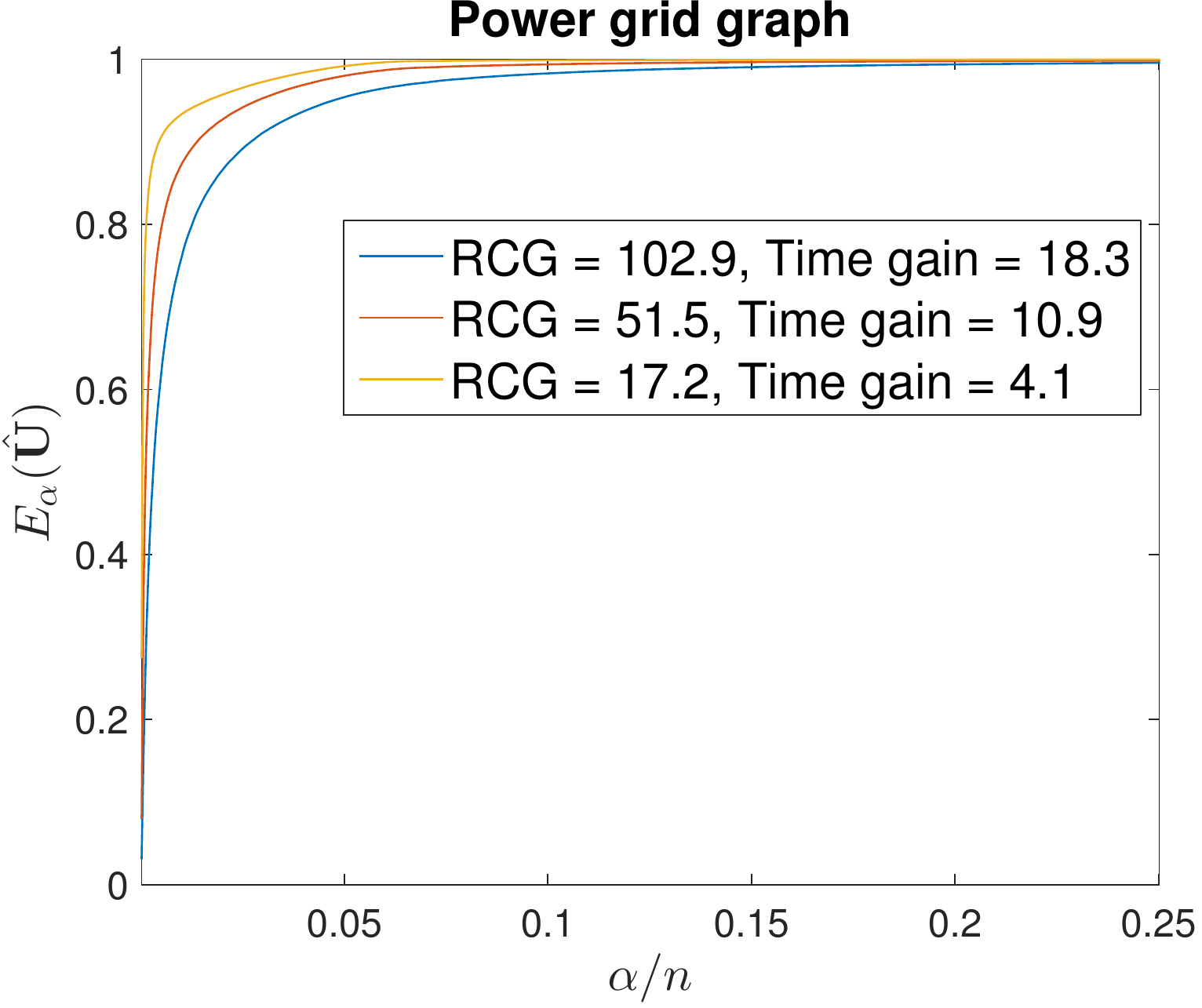}&
    
    \includegraphics[width=0.31\textwidth]{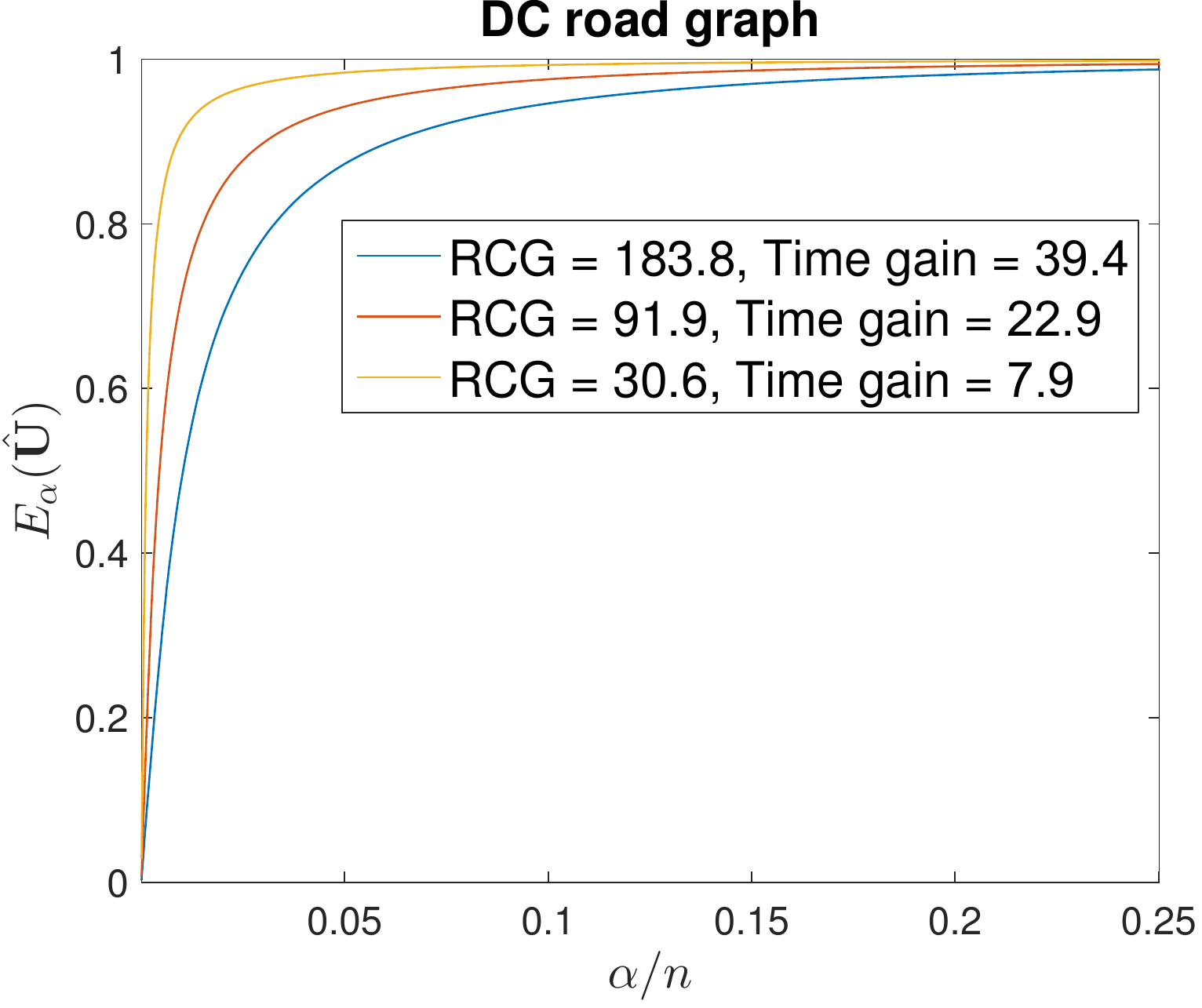}
  \end{tabular}
\caption{Relevance of several approximate \modLLM{FGFT}s computed for various graphs corresponding to sensor networks. Three approximate \modLLM{FGFT}s (corresponding to different colors) are considered for each of the three graphs. The results are shown for different values of the parameter $\alpha \in [0,\frac{n}{4}] $.}
\label{fig:DCtol}
\end{figure*}

\subsection{Experiments on real sensor graphs}

 Let us now consider an application of the approximate diagonalization method to graphs representing real sensor networks. The main idea of this experiment is to compare, for several real-world graphs corresponding to physical infrastructures: a) the approximate \modLLM{FGFT} $\hat{\mathbf{U}}_{\text{Givens}\scriptscriptstyle //}$ obtained for three different number of Givens rotations $J \in \{ n\log n, 2n\log n, 6n\log n\}$; b) to the true Fourier transform. To this end, we consider the three following graphs:
\begin{itemize}
\item {\bf Minnesota road graph:}  a classical graph 
made of $n=2642$ nodes, generated with the GSPBOX \cite{GSP}.
\item {\bf US power grid graph:} a graph representing the power grid of the western United States, made of $n=4941$ nodes. This graph was previously studied in \cite{Watts1998}, and we obtained it from \cite{PG}.
\item {\bf DC road graph:} a graph representing the road network of the District of Columbia, made of $n=9559$ nodes. This graph was previously used for a shortest path challenge \cite{DIMACS2007}, and we obtained it from \cite{DC}. \modRG{The storage of its full Fourier matrix requires 700MB.}
\end{itemize}


{\noindent \bf Performance metric.}
In order to evaluate the \modLLM{accuracy} of an approximate Fourier matrix $\hat{\mathbf{U}} = \mathbf{S}_1\dots \mathbf{S}_K$ with respect to the true Fourier matrix $\mathbf{U}$, \modLLM{we use a measure quite similar to $\text{err}_d(\hat{\mathbf{U}})$, but that is a bit more informative. Its expression is the following}:
$$
E_\alpha(\hat{\mathbf{U}}) \triangleq\frac{\big\Vert (\hat{\mathbf{U}}^T\mathbf{U})_{|i-j|\leq \alpha} \big\Vert_F^2}{\big\Vert \hat{\mathbf{U}}^T\mathbf{U} \big\Vert_F^2 } = \frac{1}{n}\left\Vert (\hat{\mathbf{U}}^T\mathbf{U})_{|i-j|\leq \alpha} \right\Vert_F^2,
$$
where $\mathbf{B}\triangleq \mathbf{A}_{|i-j|\leq \alpha}$ is a matrix of the same size as $\mathbf{A}$, with $b_{ij} = a_{ij}$ if $|j-i|\leq \alpha$, and $b_{ij}=0$ otherwise ($\mathbf{B}$ is equal to $\mathbf{A}$ in a band of width $2\alpha +1$ centered around the diagonal and null elsewhere). We argue that the quantity $E_\alpha(\hat{\mathbf{U}})$ is a good measure for the relevance of the approximate Fourier transform $\hat{\mathbf{U}}$. Indeed, it corresponds to the fraction of the total energy of the matrix $\hat{\mathbf{U}}^T\mathbf{U}$ that is contained within a band of width $2\alpha +1$ around the diagonal. Said otherwise, if we imagine a signal $\mathbf{x} \in \mathbb{R}^n$ whose true spectrum is a Dirac located at the $m$th frequency ($\mathbf{x} = \mathbf{U}\boldsymbol{\delta}_m$), $E_\alpha(\hat{\mathbf{U}})$ correspond to the fraction of the energy of the approximated spectrum $\hat{\boldsymbol{\delta}}_m \triangleq \hat{\mathbf{U}}^T\mathbf{x}$ that is within the neighborhood of width $\alpha$ of $m$ (in expectation over the choice of the frequency $m$). The measure $E_\alpha(\hat{\mathbf{U}})$ is defined for $\alpha \geq 0$, it is monotonically increasing with respect to the parameter $\alpha$, we have $0\leq E_\alpha(\hat{\mathbf{U}})\leq 1$, and the higher it is, the better is the approximate Fourier matrix $\hat{\mathbf{U}}$. Moreover, we have $E_\alpha(\mathbf{U}) = 1$, $\forall \alpha \geq 0$. In summary the quantity $E_\alpha(\hat{\mathbf{U}})$ does not measure how close the approximation $\hat{\mathbf{U}}$ is of $\mathbf{U}$ in terms of classical relative error, but  how close to it is in terms of physical interpretation of the approximate spectrum. \modLLM{Note that the $E_\alpha(\hat{\mathbf{U}})$ measure gives more information than the previously used $\text{err}_d(\hat{\mathbf{U}})$, since it describes quantitatively the energy dispersion of the approximate \modLLM{FGFT}. However, the downside of this measure is that since it depends on the parameter $\alpha$, it requires a plot to be described, whereas a single number was sufficient for $\text{err}_d(\hat{\mathbf{U}})$.}

\modLLM{Figure}~\ref{fig:DCtol} shows $E_\alpha(\hat{\mathbf{U}})$ versus $\alpha/n$, for the three considered graphs. Several comments are in order:
\begin{itemize}
\item First, it is clear that the higher the Relative Complexity Gain RCG, the lower the relevance measure $E_\alpha(\hat{\mathbf{U}})$, for any graph and any value of $\alpha$. This is quite expected and shows that there exists a compromise between computational efficiency and relevance of the approximate \modLLM{FGFT}s.
\item Second, the results seem a bit better for the power grid graph than for the two other graphs, that show very similar performance and both correspond to road networks. For example, if we consider the approximate \modLLM{FGFT}s made by the product of $J = n \log n$ Givens rotations (the blue curves), we can see that we have $E_\alpha(\hat{\mathbf{U}}) \approx 0.8$ for $\alpha/n = 0.05$ for the two road networks, whereas we have $E_\alpha(\hat{\mathbf{U}}) \approx 0.9$ for the same value of $\alpha/n$ for the power grid network. This is an interesting observation that may mean that the inherent structure of power grid networks is more suited to approximate \modLLM{FGFT}s than the structure of road networks. This hypothesis would of course require a more comprehensive study to be statistically validated.
\item Third, we can see that the larger the graph, the faster the approximate \modLLM{FGFT}s for the same error. Indeed, for the DC road graph (of size $n = 9559$), the fastest approximate \modLLM{FGFT} (the blue curve) exhibits an RCG of $183.8$ and an actual time gain of $39.4$. Its performance in terms of the relevance measure $E_\alpha(\hat{\mathbf{U}})$ is almost the same as an approximate \modLLM{FGFT} on the Minnesota road graph (of size $n=2642$) that exhibits an RCG of $10$ and an actual time gain of only $1.8$. This means that approximate \modLLM{FGFT}s of complexity $\mathcal{O}(n\log n)$ are relevant on these examples.
\end{itemize} 

\modRG{The storage of the FGFT with an RCG of 30 only requires 60MB, compared to the 700MB of the original Fourier matrix.}
\begin{figure}
\centering
\includegraphics[width=\columnwidth]{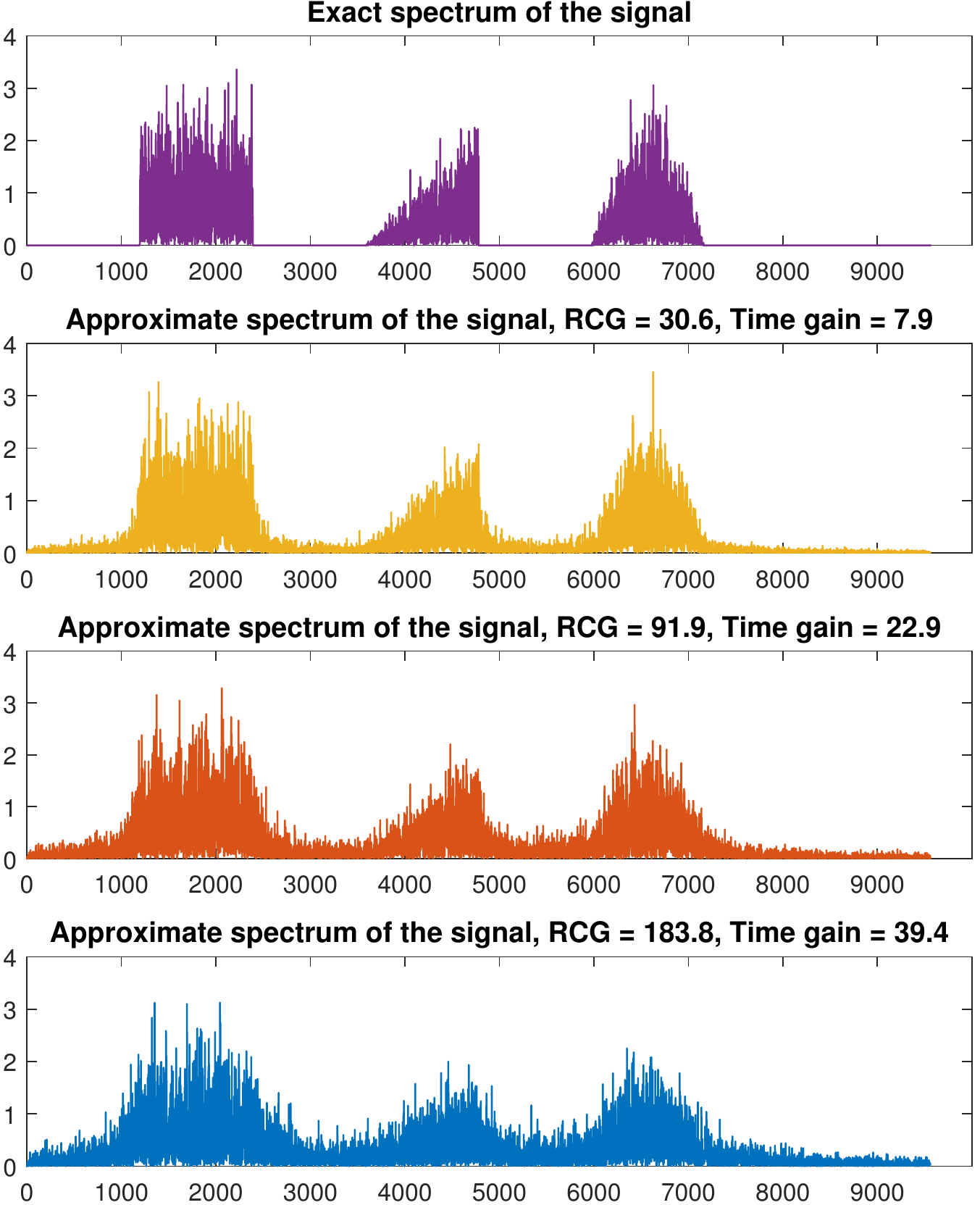}
\caption{Example of spectrum computations on the DC road graph. 
 }
\label{fig:spectrumexamplefull}
\end{figure}

Finally, let us illustrate the expected tradeoff between accuracy and computational efficiency of the $\hat{\mathbf{U}}_{\text{Givens}\scriptscriptstyle //}$ approximate \modLLM{FGFT}. To this end, we consider the DC road graph and generate a random structured spectrum $\mathbf{y} \in \mathbb{R}^{9559}$, transform it to the node domain to get $\mathbf{x} = \mathbf{Uy}$, and compute its approximate spectrum $\hat{\mathbf{y}} = \hat{\mathbf{U}}_{\text{Givens}\scriptscriptstyle //}^T\mathbf{x}$. Results are shown in \modLLM{Figure}~\ref{fig:spectrumexamplefull} for different RCG values. As expected, the larger the RCG value, the less accurate is the approximate spectrum; yet the approximation gives good results with an actual time gain of $7.9$, and fair results with an actual time gain of $22.9$.

\section{Application example: Filtering}
\label{sec:filtering}
In this last section, as an example of the versatility of approximate \modLLM{FGFTs}, we focus on one of its many possible applications: graph signal filtering. 
Given a graph signal $\mathbf{x}$, its filtered version $\mathbf{y}$ is computed as
$$
\mathbf{y} = \underbrace{\mathbf{U}\mathbf{H}\mathbf{U}^T}_{\mathbf{G}} \mathbf{x},
$$
where $\mathbf{H}$ is a diagonal matrix whose diagonal entries are the frequency response of the filter and we denote $\mathbf{G} \triangleq \mathbf{U}\mathbf{H}\mathbf{U}^T$ the filtering matrix. We see that filtering requires to apply the Fourier transform and its inverse (here, its transpose as $\mathbf{U}$ is orthogonal). An efficient approximation of the Fourier matrix would thus allow to have a complexity gain on the filtering that is of the order of the relative complexity gain RCG. \\

{\noindent \bf Filtering with an approximate \modLLM{FGFT}.}  We consider the approximate \modLLM{FGFT} $\hat{\mathbf{U}}_{\text{Givens}\scriptscriptstyle //}$
to perform filtering: the approximate filtering matrix thus reads:
$$\hat{\mathbf{G}}_{\text{\modLLM{FGFT}}} \triangleq \hat{\mathbf{U}}_{\text{Givens}\scriptscriptstyle //}\mathbf{H}\hat{\mathbf{U}}_{\text{Givens}\scriptscriptstyle //}^T.$$

 {\noindent \bf Filtering with polynomial approximations.} We compare this proposition to the usual approximate filtering method 
based on the polynomial approximation of the filter's frequency response~\cite{Hammond2011}. Starting from the frequency response of an ideal filter $h(\lambda)$ (where $\lambda$ represents the continuous frequency variable), the idea is to approximate it with a polynomial function $\hat{h}(\lambda) = \sum_{i=0}^p \alpha_i \lambda^i$, such that the approximate filtering matrix has the expression 
$$\hat{\mathbf{G}}_{\text{poly}} \triangleq \mathbf{U}\hat{h}(\boldsymbol{\Lambda})\mathbf{U}^T = \sum_{i=0}^p \alpha_i \mathbf{L}^i$$
 where $\hat{h}$ is applied entry-wise to the diagonal entries of $\boldsymbol{\Lambda}$. This provides a computationally efficient way of applying the filter, since its application cost is $p(\left\Vert \mathbf{L} \right\Vert_0 + n) + n$ arithmetic operations. With this polynomial approximation method, the Relative Complexity Gain reads $\text{RCG} = \frac{n^2}{p(\left\Vert \mathbf{L} \right\Vert_0 + n) + n}$. \modLLM{Note that obtaining the coefficients $\alpha_{i}$ characterizing $\hat{\mathbf{G}}_{\text{poly}}$ for a given filter costs $\mathcal{O}(p)$ arithmetic operations, and one needs to recompute these coefficients as soon as the filter changes. On the other hand, with the greedy diagonalizations, obtaining $\hat{\mathbf{G}}_{\text{\modLLM{FGFT}}}$ costs $\mathcal{O}(nJ\log n)$ (as denoted in Table~\ref{tab:complexity}), but can be  done once and for all and reused for as many filters as needed, as long as the same graph is used.} \\



{\noindent \bf Application to denoising.} 
We perform a first experiment on the Minnesota road network graph presented in \modLLM{Section~}\ref{sec:sensor}. 
We will compare the performance of the two approximation methods for three RCG values: 35, 17 and 12. This corresponds to a number of Givens rotations $J=50 000, 100 000, 150 000$ for the \modLLM{FGFT}-based approximation, and an order $p=14, 28, 40$ for the polynomial-based approximation. 
In this first experiment, a low-frequency signal $\mathbf{x}$ is generated randomly in the graph Fourier domain: coefficients of its spectrum $\mathbf{y}$ are independent and follow a normal distribution $y_i \sim \mathcal{N}(0,\theta_i)$ with $\theta_i = 1$ if $i\leq 300$, $\theta_i = 0$ otherwise, where $\lambda_i$ is the $i$th eigenvalue of the Laplacian. We then obtain a signal in the graph node domain $\mathbf{x} = \mathbf{U}\mathbf{y}$, that is corrupted with a white Gaussian noise $\mathbf{n}$ with $n_i \sim \mathcal{N}(0,\sigma)$, to get the noisy signal $\tilde{\mathbf{x}} = \mathbf{x} + \mathbf{n}$. The noisy signal is then low-pass filtered in order to remove part of the noise. The filter used has a frequency response given by  $h_i = 1$ if $i\leq 300$, $h_i = 0$ otherwise, where $h_i$ is the $i$th entry on the diagonal of the filter matrix $\mathbf{H}$. Filtering is done with the true Fourier matrix $\mathbf{U}$ (with the exact filter  $\mathbf{G}\tilde{\mathbf{x}}=\mathbf{U}\mathbf{H}\mathbf{U}^T\tilde{\mathbf{x}}$), to which we compare $\hat{\mathbf{G}}_{\text{poly}}\tilde{\mathbf{x}}$ and $\hat{\mathbf{G}}_{\text{\modLLM{FGFT}}}\tilde{\mathbf{x}}$.



\begin{table*}[tb]                                              
\centering                                                     
\begin{tabular}{lccccccc}                                   
\toprule                                                       
 &  $\sigma=0.2$ & $\sigma=0.25$ & $\sigma=0.3$ & $\sigma=0.4$ & $\sigma=0.5$ & $\sigma=0.6$ \\
\midrule                                                          
Noisy signal  & 4.50 & 2.61 & 0.98 & -1.45 & -3.47 & -4.99 \\                         
\hline
                                                         
Filtered with $\mathbf{U}$ & 13.97 & 12.08 & 10.40 & 7.97 & 6.00 & 4.48 \\      


\hline 

%
%
%


& $\hat{\mathbf{G}}_{\text{poly}}|\hat{\mathbf{G}}_{\text{\modLLM{FGFT}}}$
& $\hat{\mathbf{G}}_{\text{poly}}|\hat{\mathbf{G}}_{\text{\modLLM{FGFT}}}$
& $\hat{\mathbf{G}}_{\text{poly}}|\hat{\mathbf{G}}_{\text{\modLLM{FGFT}}}$
& $\hat{\mathbf{G}}_{\text{poly}}|\hat{\mathbf{G}}_{\text{\modLLM{FGFT}}}$
& $\hat{\mathbf{G}}_{\text{poly}}|\hat{\mathbf{G}}_{\text{\modLLM{FGFT}}}$
& $\hat{\mathbf{G}}_{\text{poly}}|\hat{\mathbf{G}}_{\text{\modLLM{FGFT}}}$ \\

RCG = $12$, \modLLM{$T$=2.4ms \textbar 2.9ms} &  {\bf 11.84} \textbar 11.39 & {\bf 10.70} \textbar 10.25 & {\bf 9.43} \textbar 9.05 & {\bf 7.44} \textbar 7.16 & {\bf 5.72} \textbar 5.46 & {\bf 4.34} \textbar 4.11 \\     

RCG = $17$, \modLLM{$T$=1.7ms \textbar 2.1ms}  & 0.0017  \textbar {\bf 10.35} & 8.59 \textbar  {\bf 9.43} & 7.82 \textbar {\bf 8.44} & 6.49 \textbar {\bf 6.79} & 5.16  \textbar  {\bf 5.18} &  {\bf 4.01} \textbar 3.90 \\

RCG = $35$, \modLLM{$T$=1.0ms \textbar 1.1ms}  & {\bf 8.10} \textbar 7.45 & {\bf 7.72} \textbar 7.00 & {\bf 7.10} \textbar 6.39 & {\bf 6.08} \textbar 5.25 & {\bf 4.89} \textbar 4.21 & {\bf 3.87} \textbar 3.10 \\                                       
               

\bottomrule                                                       
\end{tabular}                                                  
~\vspace{1mm}\caption{
Filtering results, the SNR in dB (mean over $100$ realizations) is given, and the best result is written in bold. \modLLM{The average filtering time $T$ (in milliseconds) for both method at each RCG is also given. 
} 
}
\label{tab:minnesota}                                 
\end{table*}

\begin{figure}[tb]
\centering
\includegraphics[width=\columnwidth]{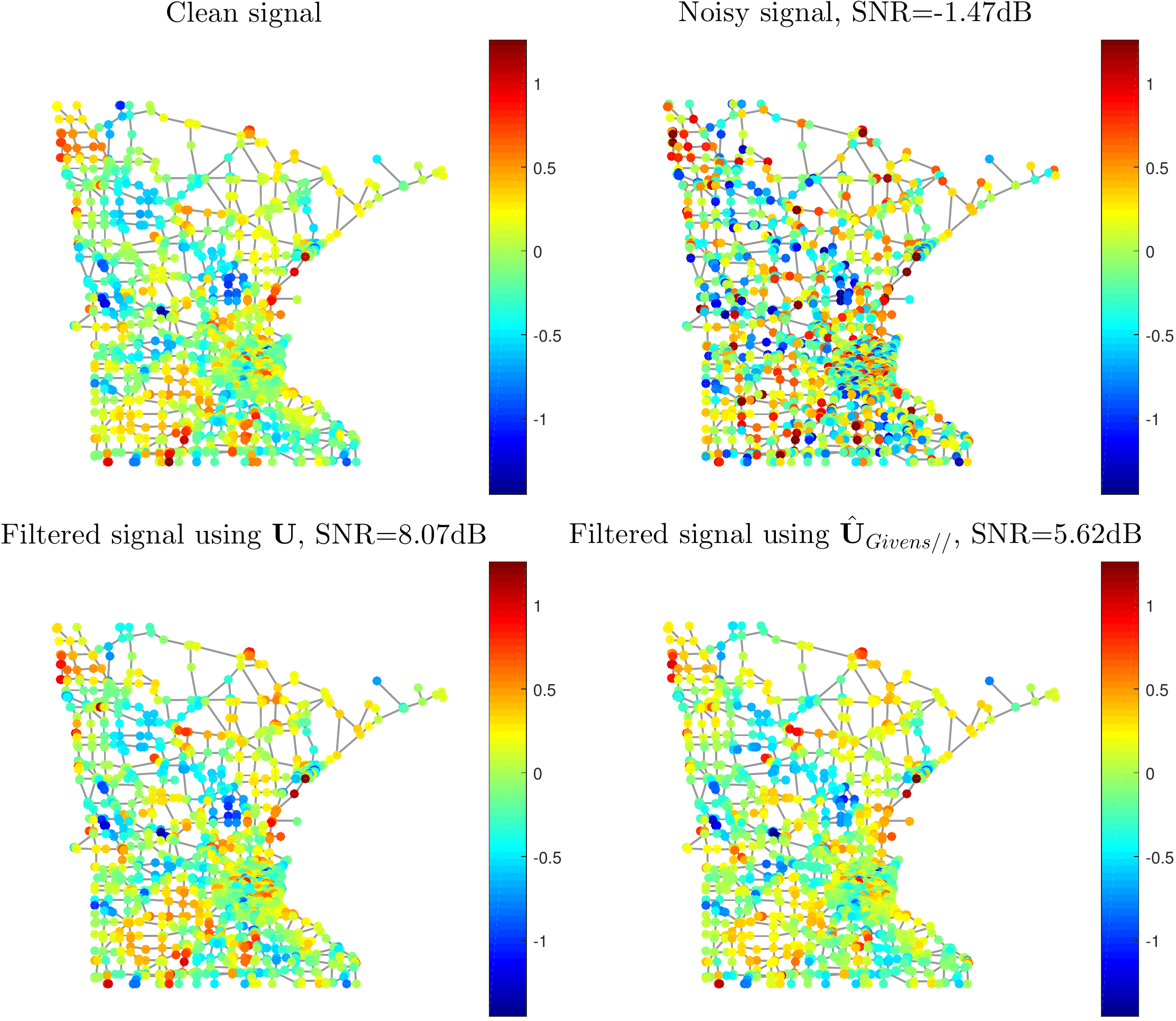}
\caption{Filtering example on the Minnesota road graph, the approximate \modLLM{FGFT} $\hat{\mathbf{U}}_{\text{Givens}\scriptscriptstyle //}$ is approximately $35$ times more efficient than the true Fourier matrix $\mathbf{U}$.}
\label{fig:denoisingexample}
\end{figure}

First of all, for illustration purposes, an example of graph filtering with  $\hat{\mathbf{U}}_{\text{Givens}\scriptscriptstyle //}$ on a realisation of noisy signal with $\sigma=0.4$ is shown in  \modLLM{Figure}~\ref{fig:denoisingexample}. Moreover, for different noise levels $\sigma$, average 
 results over $100$ realizations  are given in \modLLM{Table~}\ref{tab:minnesota}, in terms of Signal to Noise Ratio (SNR) in decibel (dB):  
$$\text{SNR} = 10\log_{10}\left(\frac{\big\Vert \hat{\mathbf{x}} \big\Vert_2^2}{\big\Vert \mathbf{x} - \hat{\mathbf{x}} \big\Vert_2^2}\right).$$
 Both methods 
 show 
 similar performance at constant RCG\modLLM{, and their run time are comparable}. 
%
%
 \textcolor{black}{However, approximate \modLLM{FGFT}s can be considered as more versatile than polynomial approximation, in the sense that they give access to an approximate spectrum, and thus allow to perform in a computationally efficient way any task requiring access to this spectrum, which is not the case of polynomial approximations that can be used for filtering only. Besides, the performance of polynomial approximations depends heavily on the considered type of filter, as is discussed in the next paragraph.}\\



{\noindent \bf Comparison of the approximate filtering matrices.}  In order to better understand the relative strengths and weaknesses of the proposed method compared to polynomial approximation, we now analyze the approximate filtering matrix $\hat{\mathbf{G}}_\text{\modLLM{FGFT}}$ 
and compare it to 
$\hat{\mathbf{G}}_\text{poly}$ obtained with equal relative complexities. 
We evaluate the accuracy of these approximate filtering matrices, generically written $\hat{\mathbf{G}}$, by comparing them with the exact filtering matrix $\mathbf{G}$, with the following measure: 
$$\frac{\big\Vert \mathbf{G} - \hat{\mathbf{G}}\big\Vert_F}{\left\Vert \mathbf{G} \right\Vert_F}.$$
In this second experiment, we consider a community graph of size $n=2048$ made of $23$ communities, and two different types of filters, representing two extreme situations with respect to polynomial approximation: an ideal low-pass filter and a \modNT{heat kernel (exponentially decreasing filter)}. The ideal low-pass filter has a frequency response equal to one for the first $r$ frequencies and zero after that (we consider that the $r$th eigenvalue is perfectly known for the polynomial approximation). We consider three choices for the index $r$ of the cutting frequency: $r=7$ which corresponds to a region where the density of graph frequencies is high, $r=23$ which corresponds to a region where there is a gap in the graph frequencies density (it corresponds to the number of communities) and $r=1000$ which corresponds to a region where the density of graph frequencies is lower. 
The \modNT{heat kernel} has a frequency response $h(\lambda) = \exp(-\lambda)$ 

\begin{table}[tb]                                                       
\centering                                                              
\begin{tabularx}{1\columnwidth}{Xccc}                                      
\toprule                                                       
&  $\text{RCG} = 35$  & $\text{RCG} = 17$ & $\text{RCG} = 12$   \\ 
& $\hat{\mathbf{G}}_{\text{poly}}|\hat{\mathbf{G}}_{\text{\modLLM{FGFT}}}$
& $\hat{\mathbf{G}}_{\text{poly}}|\hat{\mathbf{G}}_{\text{\modLLM{FGFT}}}$
& $\hat{\mathbf{G}}_{\text{poly}}|\hat{\mathbf{G}}_{\text{\modLLM{FGFT}}}$
\\
\midrule 
\scriptsize Ideal low-pass $r=\scriptstyle7$ &  {\bf 1.00}\textbar   1.03  & 1.00\textbar {\bf 0.94} & 1.39\textbar {\bf 0.93}   \\
\scriptsize Ideal low-pass $r=\scriptstyle 23$ & 1.00\textbar  {\bf 0.64}  & 0.78\textbar {\bf 0.34} & 0.42\textbar {\bf 0.23}   \\
\scriptsize Ideal low-pass $r= \scriptstyle 1000$ & 0.68\textbar  {\bf 0.25}  & 0.48\textbar {\bf 0.19} & 0.50\textbar {\bf 0.16}   \\
\scriptsize Exponential & {\bf 0.05}\textbar 0.09 & {\bf 1e-4}\textbar   0.06   &{\bf 0.00}\textbar 0.05  \\
\bottomrule                                                             
\end{tabularx}                                                           
~\vspace{0mm}\caption{Comparison \modLLM{(relative Frobenius norm error)} between the filtering operators obtained with the approximate \modLLM{FGFT} and with polynomial approximations, the best result being written in bold. 
} 
\label{tab:compcheby}                                              
\end{table}

Table~\ref{tab:compcheby} shows the resulting relative approximation errors for various RCG values. For the ideal low-pass filter with $r=7$
, the two approximate filtering methods 
perform quite poorly, with a small advantage for the \modLLM{FGFT}-based approximation. 
This can be explained by the fact that this filtering operation is intrinsically difficult to approximate because of the high density of eigenvalues around the cutting frequency. On the other hand, with $r=23$, results are better for the \modLLM{FGFT}-based approximation, but remain quite poor for the polynomial approximations. Finally, for $r=1000$
, results are better for both methods while the \modLLM{FGFT}-based approximation shows a clear advantage for all tested values of RCG. Conversely, for the exponential filter, polynomial approximations give much better results than approximate \modLLM{FGFT}s (it is essentially perfect). This can be explained by the fact that the performance of polynomial approximations depends heavily on the type of filter: it works well for filters whose frequency response is well approximated by polynomials of low degree (typically the exponential filter), and works not as well for filters whose frequency response is not easily approximated (typically the ideal low-pass filter). \modLLM{Indeed, polynomials of low degree can approximate well only relatively smooth functions \cite{Dupont1980}}. On the other hand, the performance of approximate \modLLM{FGFT}s seems less dependent on the type of filters. 

It is thus expected that in general the polynomial filter approximation is the best method for filtering signals on graphs, but in specific situations (such as the ideal low-pass filter), it may be outperformed by an \modLLM{FGFT}-based approximation. Interestingly, the ideal low-pass filter is important in certain applications such as spectral clustering \cite{Tremblay2016}, where one needs to apply an ideal low-pass filter with $r$ being equal to the desired number of clusters. 

\section{Conclusion}
In this paper, we proposed a method to obtain approximate Fast Fourier Transforms on graphs via an approximate diagonalization of the Laplacian matrix with an approximate eigenvector matrix constrained to be a product of a relatively low number of Givens rotations. The method was first described and analyzed \modLLM{in detail}. It was then compared to another method we previously proposed on various popular graphs, and applied to sensor networks, showing promising results. Finally, an application to filtering was proposed, showing that approximate \modLLM{FGFT}s can perform as good as polynomial filter approximations for certain filters.

In the future, it would be very interesting to reduce the gap between actual time gain and RCG, to have approximate \modLLM{FGFT}s whose benefits are closer to the theoretical gain, whose usability would be enhanced.

It would also be beneficial to the method to imagine clever ways to update the approximate \modLLM{FGFT} when the graph changes, without requiring to recompute a complete approximate diagonalization. This could indeed open up the method to new applications (other than sensor networks) where the graphs are varying by nature (social networks, movie ratings, etc.).

We could also imagine applying the method to approximately diagonalize other symmetric matrices besides graph Laplacians. For example, applying it directly to a graph filtering matrix $\mathbf{G}$ (computed for example by polynomial approximation, as defined in \modLLM{Section~}\ref{sec:filtering}, \textcolor{black}{to avoid any costly diagonalization}) could allow for a faster application of the filter, which is of interest if the filter has to be applied a great number of times.

Furthermore, an efficient implementation that would reduce the approximate diagonalization time of the proposed methods would be very interesting. Indeed, efficient numerical linear algebra routines exist that implement the Jacobi eigenvalue algorithm, so that modifying them directly could potentially greatly reduce this factorization time. \textcolor{black}{Moreover, the use of more advanced techniques similar to those introduced in \cite{Teneva2016} for the pMMF method may also be used to accelerate the diagonalization.} This way, the cost of approximate diagonalization (whose theoretical time complexity is $\mathcal{O}(n^2\log n)$ for $J = \mathcal{O}(n\log n)$ Givens rotations) may even become competitive with that of the cost of applying the true Fourier transform (whose theoretical time complexity is $\mathcal{O}(n^2)$). This would open new fields of applications for the proposed methods.

\modLLM{Finally, a theoretical analysis of the proposed greedy diagonalization would be of great interest. For example, a result linking a graph structure and the minimum number of Givens rotations required to diagonalize its Laplacian up to a prescribed accuracy would be very useful to get a deeper understanding of the method. More generally, characterizing which symmetric matrices admit good approximate diagonalization by products of few Givens (or more generally few sparse) factors is an interesting question for approximation theory.}


%

\appendices

\section{Efficient greedy strategy}
\label{app:efficientdiago}
We give here a computationally efficient version of the algorithm of \modLLM{Figure}~\ref{algo:givens_sum}. 
\begin{figure}[htbp]
\centering 
\scriptsize
\begin{boxedalgorithmic}[1] 
\scriptsize
\REQUIRE{Matrix $\mathbf{L}$, desired number of factors $J$.}\STATE $\mathbf{L}_1 \leftarrow \mathbf{L}$
\FOR{$r=1$ to $n$}
\FOR{$s=r+1$ to $n$}
\STATE $c_{rs} \leftarrow  -|l_{rs}^j|$
\ENDFOR
\STATE $d(r) \leftarrow \min \mathbf{C}(r,:)$
\STATE $e(r) \leftarrow \text{argmin} \: \mathbf{C}(r,:)$ 
\ENDFOR
\FOR{$j=1$ to $J$}
\STATE $p \leftarrow \text{argmin} \: d(r)$
\STATE $q \leftarrow e(p)$ 
\STATE $\theta \leftarrow \frac{1}{2}\arctan(\frac{l^j_{qq} - l^j_{pp}}{2l^j_{pq}}) + \frac{\pi}{4}$
\STATE  $\mathbf{S}_j \leftarrow \mathbf{G}_{p,q,\theta} $
\STATE  $\mathbf{L}_{j+1} \leftarrow \mathbf{S}_j^T\mathbf{L}_j\mathbf{S}_j$
\FOR{$r=p,q$}
\FOR{$s=r+1$ to $n$}
\STATE $c_{rs} \leftarrow  -|l_{rs}^j|$
\ENDFOR
\STATE $d(r) \leftarrow \min \mathbf{C}(r,:)$
\STATE $e(r) \leftarrow \text{argmin} \: \mathbf{C}(r,:)$ 
\ENDFOR
\FOR{$s=p,q$}
\FOR{$r=1$ to $s-1$}
\STATE $c_{rs} \leftarrow  -|l_{rs}^j|$
\IF{$c_{rs}<d(r)$}
\STATE $d(r)\leftarrow c_{rs}$
\STATE $e(r)\leftarrow s$
\ELSE
\IF{$e(r)= s$}
\STATE $d(r) \leftarrow \min \mathbf{C}(r,:)$
\STATE $e(r) \leftarrow \text{argmin} \: \mathbf{C}(r,:)$ 
\ENDIF
\ENDIF
\ENDFOR
\ENDFOR
\ENDFOR
\STATE  $\hat{\boldsymbol{\Lambda}} \leftarrow \text{diag}(\mathbf{L}_{J+1})$
\modfinal{\STATE Sort diagonal entries of $\hat{\boldsymbol{\Lambda}}$ in increasing order. Reorder columns of $\mathbf{S}_J$ accordingly.}
\ENSURE Sparse orthogonal factors  $\{\mathbf{S}_j\}_{j=1}^{J}$; diagonal factor $\hat{\boldsymbol{\Lambda}}$.
\end{boxedalgorithmic}
\label{algo:givens_full_efficient}
\caption{{\bf Truncated Jacobi algorithm (efficient version) }: Algorithm for the approximate greedy diagonalization of the Laplacian matrix (efficient version).}
\end{figure}

\section*{Acknowledgment}

The authors would like to thank Igal Rozenberg for discussions and collaboration on the topic of graph signal processing during his internship with the team.

\ifCLASSOPTIONcaptionsoff
  \newpage
\fi



\bibliographystyle{IEEEtran}
\bibliography{IEEEabrv,biblio}

\begin{thebibliography}{10}
\providecommand{\url}[1]{#1}
\csname url@samestyle\endcsname
\providecommand{\newblock}{\relax}
\providecommand{\bibinfo}[2]{#2}
\providecommand{\BIBentrySTDinterwordspacing}{\spaceskip=0pt\relax}
\providecommand{\BIBentryALTinterwordstretchfactor}{4}
\providecommand{\BIBentryALTinterwordspacing}{\spaceskip=\fontdimen2\font plus
\BIBentryALTinterwordstretchfactor\fontdimen3\font minus
  \fontdimen4\font\relax}
\providecommand{\BIBforeignlanguage}[2]{{%
\expandafter\ifx\csname l@#1\endcsname\relax
\typeout{** WARNING: IEEEtran.bst: No hyphenation pattern has been}%
\typeout{** loaded for the language `#1'. Using the pattern for}%
\typeout{** the default language instead.}%
\else
\language=\csname l@#1\endcsname
\fi
#2}}
\providecommand{\BIBdecl}{\relax}
\BIBdecl

\bibitem{Bondy1976}
J.~A. Bondy and U.~S.~R. Murty, \emph{Graph theory with applications}.\hskip
  1em plus 0.5em minus 0.4em\relax Elsevier Science Ltd/North-Holland, 1976,
  vol. 290.

\bibitem{Shuman2013}
D.~I. Shuman, S.~K. Narang, P.~Frossard, A.~Ortega, and P.~Vandergheynst, ``The
  emerging field of signal processing on graphs: Extending high-dimensional
  data analysis to networks and other irregular domains,'' \emph{Signal
  Processing Magazine, IEEE}, vol.~30, no.~3, pp. 83--98, 2013.

\bibitem{Sandryhaila2013}
A.~Sandryhaila and J.~Moura, ``Discrete signal processing on graphs,''
  \emph{Signal Processing, IEEE Transactions on}, vol.~61, no.~7, pp.
  1644--656, 2013.

\bibitem{CooleyTukey1965}
J.~Cooley and J.~Tukey, ``An algorithm for the machine calculation of complex
  {F}ourier series,'' \emph{Mathematics of Computation}, vol.~19, no.~90, pp.
  297--301, 1965.

\bibitem{Morgenstern1975}
\BIBentryALTinterwordspacing
J.~Morgenstern, ``The linear complexity of computation,'' \emph{J. ACM},
  vol.~22, no.~2, pp. 184--194, Apr. 1975. [Online]. Available:
  \url{http://doi.acm.org/10.1145/321879.321881}
\BIBentrySTDinterwordspacing

\bibitem{Lemagoarou2016a}
L.~Le~Magoarou and R.~Gribonval, ``Flexible multilayer sparse approximations of
  matrices and applications,'' \emph{IEEE Journal of Selected Topics in Signal
  Processing}, vol.~10, no.~4, pp. 688--700, 2016.

\bibitem{Lemagoarou2016}
------, ``Are there approximate fast fourier transforms on graphs?'' in
  \emph{2016 IEEE International Conference on Acoustics, Speech and Signal
  Processing (ICASSP)}, March 2016, pp. 4811--4815.

\bibitem{Givens1958}
\BIBentryALTinterwordspacing
W.~Givens, ``\BIBforeignlanguage{English}{Computation of plane unitary
  rotations transforming a general matrix to triangular form},''
  \emph{\BIBforeignlanguage{English}{Journal of the Society for Industrial and
  Applied Mathematics}}, vol.~6, no.~1, pp. pp. 26--50, 1958. [Online].
  Available: \url{http://www.jstor.org/stable/2098861}
\BIBentrySTDinterwordspacing

\bibitem{Jacobi1846}
C.~G.~J. Jacobi, ``\"{U}ber ein leichtes verfahren, die in der theorie der
  s\"{a}kularst\"{o}rungen vorkommenden gleichungen numerisch aufzul\"{o}sen,''
  \emph{J. reine angew. Math.}, vol.~30, pp. 51--94, 1846.

\bibitem{Rubinstein2010a}
R.~Rubinstein, M.~Zibulevsky, and M.~Elad, ``Double sparsity: Learning sparse
  dictionaries for sparse signal approximation,'' \emph{Signal Processing, IEEE
  Transactions on}, vol.~58, no.~3, pp. 1553--1564, March 2010.

\bibitem{Chabiron2014}
\BIBentryALTinterwordspacing
O.~Chabiron, F.~Malgouyres, J.-Y. Tourneret, and N.~Dobigeon,
  ``\BIBforeignlanguage{English}{Toward fast transform learning},''
  \emph{\BIBforeignlanguage{English}{International Journal of Computer
  Vision}}, pp. 1--22, 2014. [Online]. Available:
  \url{http://dx.doi.org/10.1007/s11263-014-0771-z}
\BIBentrySTDinterwordspacing

\bibitem{Sulam2016}
J.~Sulam, B.~Ophir, M.~Zibulevsky, and M.~Elad, ``Trainlets: Dictionary
  learning in high dimensions,'' \emph{IEEE Transactions on Signal Processing},
  2016.

\bibitem{Chabiron2016}
\BIBentryALTinterwordspacing
O.~Chabiron, F.~Malgouyres, H.~Wendt, and J.-Y. Tourneret, ``{Optimization of a
  Fast Transform Structured as a Convolutional Tree},'' Jan. 2016, working
  paper or preprint. [Online]. Available:
  \url{https://hal.archives-ouvertes.fr/hal-01258514}
\BIBentrySTDinterwordspacing

\bibitem{Lyu2013}
S.~Lyu and X.~Wang, ``On algorithms for sparse multi-factor {NMF},'' in
  \emph{Advances in Neural Information Processing Systems 26}, 2013, pp.
  602--610.

\bibitem{Beylkin1991}
G.~Beylkin, R.~Coifman, and V.~Rokhlin, ``Fast wavelet transforms and numerical
  algorithms i,'' \emph{Communications on pure and applied mathematics},
  vol.~44, no.~2, pp. 141--183, 1991.

\bibitem{Rokhlin1985}
\BIBentryALTinterwordspacing
V.~Rokhlin, ``Rapid solution of integral equations of classical potential
  theory,'' \emph{Journal of Computational Physics}, vol.~60, no.~2, pp. 187 --
  207, 1985. [Online]. Available:
  \url{http://www.sciencedirect.com/science/article/pii/0021999185900026}
\BIBentrySTDinterwordspacing

\bibitem{Hackbusch1999}
\BIBentryALTinterwordspacing
W.~Hackbusch, ``{A Sparse Matrix Arithmetic Based on H-matrices. Part I:
  Introduction to H-matrices},'' \emph{Computing}, vol.~62, no.~2, pp. 89--108,
  May 1999. [Online]. Available: \url{http://dx.doi.org/10.1007/s006070050015}
\BIBentrySTDinterwordspacing

\bibitem{Candes2007}
\BIBentryALTinterwordspacing
E.~Cand\`es, L.~Demanet, and L.~Ying, ``Fast computation of fourier integral
  operators,'' \emph{SIAM Journal on Scientific Computing}, vol.~29, no.~6, pp.
  2464--2493, 2007. [Online]. Available:
  \url{http://dx.doi.org/10.1137/060671139}
\BIBentrySTDinterwordspacing

\bibitem{Rusu2016}
C.~Rusu, N.~Gonzalez-Prelcic, and R.~W.~J. Heath, ``Fast orthonormal
  sparsifying transforms based on householder reflectors,'' \emph{Submitted to
  IEEE Transactions on Signal Processing}, 2016.

\bibitem{Householder1958}
\BIBentryALTinterwordspacing
A.~S. Householder, ``Unitary triangularization of a nonsymmetric matrix,''
  \emph{J. ACM}, vol.~5, no.~4, pp. 339--342, Oct. 1958. [Online]. Available:
  \url{http://doi.acm.org/10.1145/320941.320947}
\BIBentrySTDinterwordspacing

\bibitem{Lee2008}
A.~B. Lee, B.~Nadler, and L.~Wasserman, ``Treelets - an adaptive multi-scale
  basis for sparse unordered data,'' \emph{The Annals of Applied Statistics},
  vol.~2, no.~2, pp. 435--471, July 2008.

\bibitem{Cao2011}
G.~Cao, L.~Bachega, and C.~Bouman, ``The sparse matrix transform for covariance
  estimation and analysis of high dimensional signals,'' \emph{Image
  Processing, IEEE Transactions on}, vol.~20, no.~3, pp. 625--640, 2011.

\bibitem{Chen2011}
H.~Chen, S.~Zhu, and B.~Zeng, ``Design of non-separable transforms for
  directional 2-d sources,'' in \emph{Image Processing (ICIP), 2011 18th IEEE
  International Conference on}.\hskip 1em plus 0.5em minus 0.4em\relax IEEE,
  2011, pp. 3697--3700.

\bibitem{Chen2012}
H.~Chen and B.~Zeng, ``Design of low-complexity, non-separable 2-d transforms
  based on butterfly structures,'' in \emph{Circuits and Systems (ISCAS), 2012
  IEEE International Symposium on}.\hskip 1em plus 0.5em minus 0.4em\relax
  IEEE, 2012, pp. 2921--2924.

\bibitem{Kondor2014}
R.~Kondor, N.~Teneva, and V.~Garg, ``Multiresolution matrix factorization,'' in
  \emph{Proceedings of the 31st International Conference on Machine Learning
  (ICML-14)}, 2014, pp. 1620--1628.

\bibitem{Teneva2016}
N.~Teneva, P.~K. Mudrakarta, and R.~Kondor, ``Multiresolution matrix
  compression,'' in \emph{Proceedings of the 19th International Conference on
  Artificial Intelligence and Statistics}, 2016, pp. 1441--1449.

\bibitem{Sandryhaila2014}
A.~Sandryhaila and J.~M. Moura, ``Big data analysis with signal processing on
  graphs: Representation and processing of massive data sets with irregular
  structure,'' \emph{Signal Processing Magazine, IEEE}, vol.~31, no.~5, pp.
  80--90, Sept 2014.

\bibitem{Golub2000}
G.~H. Golub and H.~A. Van~der Vorst, ``Eigenvalue computation in the 20th
  century,'' \emph{Journal of Computational and Applied Mathematics}, vol. 123,
  no.~1, pp. 35--65, 2000.

\bibitem{Golub2012}
G.~H. Golub and C.~F. Van~Loan, \emph{Matrix computations}.\hskip 1em plus
  0.5em minus 0.4em\relax JHU Press, 2012, vol.~3.

\bibitem{Strassen1969}
V.~Strassen, ``Gaussian elimination is not optimal,'' \emph{Numerische
  Mathematik}, vol.~13, no.~4, pp. 354--356, 1969.

\bibitem{Coppersmith1990}
D.~Coppersmith and S.~Winograd, ``Matrix multiplication via arithmetic
  progressions,'' \emph{Journal of symbolic computation}, vol.~9, no.~3, pp.
  251--280, 1990.

\bibitem{Legall2014}
\BIBentryALTinterwordspacing
F.~L. Gall, ``Powers of tensors and fast matrix multiplication,'' in
  \emph{International Symposium on Symbolic and Algebraic Computation, {ISSAC}
  '14, Kobe, Japan, July 23-25, 2014}, 2014, pp. 296--303. [Online]. Available:
  \url{http://doi.acm.org/10.1145/2608628.2608664}
\BIBentrySTDinterwordspacing

\bibitem{Brent1985}
R.~P. Brent and F.~T. Luk, ``The solution of singular-value and symmetric
  eigenvalue problems on multiprocessor arrays,'' \emph{SIAM Journal on
  Scientific and Statistical Computing}, vol.~6, no.~1, pp. 69--84, 1985.

\bibitem{GSP}
N.~{Perraudin}, J.~{Paratte}, D.~{Shuman}, V.~{Kalofolias}, P.~{Vandergheynst},
  and D.~K. {Hammond}, ``{GSPBOX: A toolbox for signal processing on graphs},''
  \emph{ArXiv e-prints}, Aug. 2014.

\bibitem{Bolte2014}
\BIBentryALTinterwordspacing
J.~Bolte, S.~Sabach, and M.~Teboulle, ``\BIBforeignlanguage{English}{Proximal
  alternating linearized minimization for nonconvex and nonsmooth problems},''
  \emph{\BIBforeignlanguage{English}{Mathematical Programming}}, vol. 146, no.
  1-2, pp. 459--494, 2014. [Online]. Available:
  \url{http://dx.doi.org/10.1007/s10107-013-0701-9}
\BIBentrySTDinterwordspacing

\bibitem{Tremblay2016}
N.~Tremblay, G.~Puy, R.~Gribonval, and P.~Vandergheynst, ``Compressive spectral
  clustering,'' in \emph{Machine Learning, Proceedings of the Thirty-third
  International Conference (ICML 2016), June}, 2016, pp. 20--22.

\bibitem{decelle_asymptotic_2011}
A.~Decelle, F.~Krzakala, C.~Moore, and L.~Zdeborová, ``Asymptotic analysis of
  the stochastic block model for modular networks and its algorithmic
  applications,'' \emph{Phys. Rev. E}, vol.~84, no.~6, p. 066106, Dec. 2011.

\bibitem{Sandryhaila2014a}
A.~Sandryhaila and J.~M. Moura, ``Discrete signal processing on graphs:
  Frequency analysis,'' \emph{Signal Processing, IEEE Transactions on},
  vol.~62, no.~12, pp. 3042--3054, 2014.

\bibitem{Watts1998}
D.~J. Watts and S.~H. Strogatz, ``Collective dynamics of 'small-world'
  networks,'' \emph{Nature}, vol. 393, pp. 440--442, 1998.

\bibitem{PG}
\BIBentryALTinterwordspacing
``Us power grid network dataset.'' [Online]. Available:
  \url{http://konect.uni-koblenz.de/networks/opsahl-powergrid}
\BIBentrySTDinterwordspacing

\bibitem{DIMACS2007}
C.~Demetrescu, A.~Goldberg, and D.~Johnson, ``9th dimacs challenge on shortest
  paths,'' 2007, http://www.dis.uniroma1.it/challenge9/.

\bibitem{DC}
\BIBentryALTinterwordspacing
``Dc road network dataset.'' [Online]. Available:
  \url{http://www.dis.uniroma1.it/challenge9/data/tiger/}
\BIBentrySTDinterwordspacing

\bibitem{Hammond2011}
\BIBentryALTinterwordspacing
D.~K. Hammond, P.~Vandergheynst, and R.~Gribonval, ``Wavelets on graphs via
  spectral graph theory,'' \emph{Applied and Computational Harmonic Analysis},
  vol.~30, no.~2, pp. 129 -- 150, 2011. [Online]. Available:
  \url{http://www.sciencedirect.com/science/article/pii/S1063520310000552}
\BIBentrySTDinterwordspacing

\bibitem{Dupont1980}
T.~Dupont and R.~Scott, ``Polynomial approximation of functions in sobolev
  spaces,'' \emph{Mathematics of Computation}, vol.~34, no. 150, pp. 441--463,
  1980.

\end{thebibliography}
%
%
%

%

\begin{IEEEbiography}[{\includegraphics[width=1in,height=1.25in,clip,keepaspectratio]{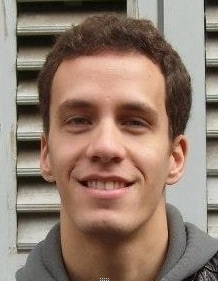}}]{Luc Le Magoarou} is a postdoctoral researcher at b\scalebox{0.8}{$<>$}com (Rennes, France). He received the Ph.D. in signal processing and the M.Sc. in electrical engineering, both from the National Institute of Applied Sciences (INSA Rennes, France), in 2016 and 2013 respectively. During his Ph.D., he was with Inria (Rennes, France), in the PANAMA research group. His main research interests lie in signal processing and machine learning, with an emphasis on computationally efficient methods and matrix factorization.
\end{IEEEbiography}
%
\begin{IEEEbiography}[{\includegraphics[width=1in,height=1.25in,clip,keepaspectratio]{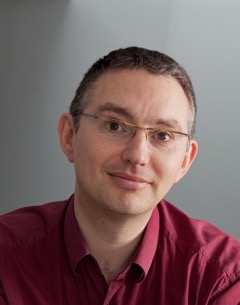}}]{R\'emi Gribonval}(FM'14)  is a Senior Researcher with Inria (Rennes, France), and the scientific leader of the PANAMA research group on sparse audio processing. A former student at  {\'E}cole Normale Sup{\'e}rieure (Paris, France), he received the Ph.D. degree in applied mathematics from Universit{\'e} de Paris-IX Dauphine (Paris, France) in 1999, and his Habilitation {\`a} Diriger des Recherches in applied mathematics from Universit{\'e} de Rennes~I (Rennes, France) in 2007. His research focuses on mathematical signal processing, machine learning, approximation theory and statistics, with an emphasis on sparse approximation, audio source separation, dictionary learning and compressed sensing. 
\end{IEEEbiography}

\begin{IEEEbiography}[{\includegraphics[width=1in,height=1.25in,clip,keepaspectratio]{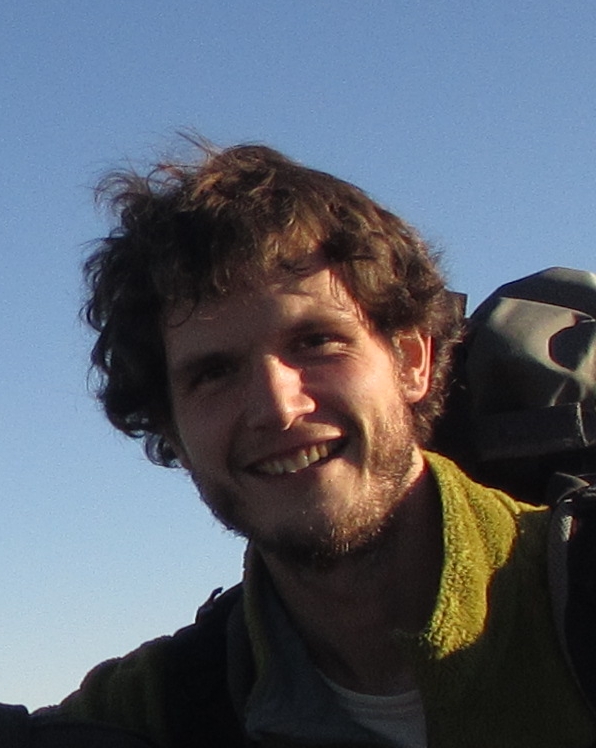}}]{Nicolas Tremblay} received his M.Sc. degree in theoretical physics with a minor in complex systems in 2009, and his PhD in physics and signal processing in 2014, from the Ecole Normale Supérieure de Lyon (ENS). He worked in 2015-2016 on a collaboration project between Inria (Rennes, France) and EPFL (Lausanne, Switzerland). He is a CNRS research fellow since 2016. His research interests are at the interface between signal processing and complex networks : graph signal processing, multiscale community detection, fast spectral algorithms for classification, and their applications to numerous networks from different fields of research, such as biological, sensor, or social networks.
\end{IEEEbiography}




\end{document}